\documentclass[12pt,a4paper]{amsart}

\author{Nils Matthes}

\date{}
\address{Max-Planck-Institut f\"ur Mathematik, Vivatsgasse 7, D-53111, Bonn, Germany}
\email{nilsmath@mpim-bonn.mpg.de}

\usepackage{geometry}
\usepackage{mathtools}
\mathtoolsset{showonlyrefs}

\usepackage{diagbox}
\usepackage{enumerate}
\usepackage{soul}

\setul{}{1pt}

\makeindex

\usepackage[T1]{fontenc}
\usepackage[english]{babel}
\usepackage[latin1]{inputenc}
\usepackage{lmodern}
\usepackage{ngerman}
\usepackage{amsmath}
\usepackage{amssymb}
\usepackage{amscd}
\usepackage{enumitem}
\usepackage[bookmarksopen]{hyperref}
\usepackage{amsthm}
\usepackage{dsfont}
\usepackage{fancyhdr}
\usepackage{mathrsfs}
\usepackage{tikz}
\usetikzlibrary{matrix,arrows}
\usepackage{pstricks}
\usepackage{graphicx}
\usepackage{color}
\usepackage{ifpdf}
\usepackage{wrapfig}
\usepackage{shuffle}

\usetikzlibrary{backgrounds}

\theoremstyle{definition}
\newtheorem{dfn}{Definition}[section]
\newtheorem{rmk}[dfn]{Remark}
\newtheorem{exmp}[dfn]{Example}

\theoremstyle{plain}
\newtheorem{prop}[dfn]{Proposition}

\newtheorem{thm}[dfn]{Theorem}
\newtheorem{lem}[dfn]{Lemma}
\newtheorem{cor}[dfn]{Corollary}

\newtheorem*{intthm}{Theorem}

\newenvironment{prf}{\begin{proof}[{\it Proof: \nopunct}]}{\end{proof}}

\def\bQ          {\mathbb Q}

\def\bC          {\mathbb C}
\def\bZ          {\mathbb Z}

\def\cO        {\mathcal O}
\def\cA           {\mathcal A}
\def\cB           {\mathcal B}
\def\cC           {\mathcal C}
\def\cF         {{\mathcal F}}

\def\cB         {{\mathcal B}}

\def\cI         { {\mathcal I}}

\def\dd       {\mathrm{d}}

\DeclareMathOperator{\SL}{SL}
\DeclareMathOperator{\im}{Im}

\DeclareMathOperator{\spn}{Span}

\numberwithin{equation}{section}
\DeclareMathOperator{\Frac}{Frac}

\newif\ifnote 
\notetrue

\allowdisplaybreaks

\def\fH          {\mathfrak{H}}

\RequirePackage{verbatim}

\makeatletter
\newwrite\bibinl@out
\newenvironment{bibtex}[1][\jobname]{%
 \immediate\openout\bibinl@out #1.bib
 \immediate\write\bibinl@out{\@percentchar generated from `\jobname' starting line \the
\inputlineno^^J}%
 \def\verbatim@processline{\immediate\write\bibinl@out{\the\verbatim@line}}%
 \@bsphack\let\do\@makeother\dospecials\catcode`\^^M\active\verbatim@start
}%
{\immediate\closeout\bibinl@out\@esphack}
\makeatother

\begin{document}
\title[Iterated integrals of quasimodular forms]{On the algebraic structure of iterated integrals of quasimodular forms}

\begin{abstract}
We study the algebra $\cI^{QM}$ of iterated integrals of quasimodular forms for $\SL_2(\bZ)$, which is the smallest extension of the algebra $QM_{\ast}$ of quasimodular forms, which is closed under integration. We prove that $\cI^{QM}$ is a polynomial algebra in infinitely many variables, given by Lyndon words on certain monomials in Eisenstein series.
We also prove an analogous result for the $M_{\ast}$-subalgebra $\cI^{M}$ of $\cI^{QM}$ of iterated integrals of modular forms.
\end{abstract}

\maketitle
\section{Introduction}
Quasimodular forms are generalizations of modular forms, which have first been introduced in \cite{KanekoZagier}, in a context motivated by mathematical physics. The $\bC$-algebra $QM_{\ast}$ of quasimodular forms for the full modular group $\SL_2(\bZ)$ can be defined, in a slightly ad hoc fashion, as the polynomial ring $\bC[E_2,E_4,E_6]$, where $E_{2k}$ denotes the normalized Eisenstein series of weight $2k$:
\begin{equation}
E_{2k}(\tau)=1-\frac{4k}{B_{2k}}\sum_{n=1}^{\infty}n^{2k-1}\frac{q^n}{1-q^n}, \quad q=e^{2\pi i\tau},
\end{equation}
where $B_{2k}$ are the Bernoulli numbers. In particular, $QM_{\ast}$ contains the algebra of modular forms $M_{\ast} \cong \bC[E_4,E_6]$.

The derivative of a quasimodular form (of weight $k$) is again a quasimodular form (of weight $k+2$); this was essentially already known to Ramanujan (cf. \cite{123}, Proposition 15). On the other hand, the integral of a quasimodular form is in general not quasimodular. For example, a primitive of $E_2$ would have to be of weight zero, but every quasimodular form of weight zero is constant.

The goal of this paper is to study the smallest algebra extension of $QM_{\ast}$, which is closed under integration. For this, the idea is to iteratively adjoin primitives to $QM_{\ast}$, which eventually leads to adjoining all (indefinite) \textit{iterated integrals}
\begin{equation} \label{eqn:iterEich}
I(f_1,\ldots,f_n;\tau)=(2\pi i)^n\idotsint\limits_{\tau \leq \tau_1 \leq \ldots \leq \tau_n\leq i\infty}f_1(\tau_1)\ldots f_n(\tau_n)\dd\tau_1\ldots\dd\tau_n,
\end{equation}
where $f_1,\ldots,f_n$ are quasimodular forms (a precise definition will be given in Section \ref{ssec:2.2}). The integrals \eqref{eqn:iterEich} have first been studied by Manin \cite{Man} and later by Brown \cite{Brown:MMV} and Hain \cite{Hain:HodgeDeRham}, in the case where all the $f_i$ are modular forms.\footnote{More precisely, Manin only defined iterated integrals of cusp forms, and the extension to all modular forms is due to Brown.} In all of these treatments, the focus lies rather on arithmetic aspects of these iterated integrals, for example their special values at cusps of the upper half-plane. By contrast, we study them solely as holomorphic functions of $\tau$. It is also worth noting that even in the modular case, the iterated integrals we study in the present paper are slightly more general than the ones introduced in \cite{Brown:MMV,Hain:HodgeDeRham,Man}. For example, if $f(\tau)$ is a modular form of weight $k$, then the integral $\int_{\tau}^{i\infty}f(\tau_1)\tau_1^n\dd\tau_1$ is an iterated integral of modular forms in the sense of the present paper \textit{for every} $n \geq 0$, while \cite{Brown:MMV,Hain:HodgeDeRham,Man} also require $n \leq k-2$.

Now let $\cI^{QM}$ be the $QM_{\ast}$-algebra generated by all the integrals \eqref{eqn:iterEich}, which is the smallest algebra extension of $QM_{\ast}$, closed under integration. It turns out that $\cI^{QM}$ is not finitely generated, but still has a manageable structure, which is captured by the notion of shuffle algebra (which is just the graded dual of the tensor algebra with a certain commutative multiplication, the so-called shuffle product) \cite{Reu}.
More precisely, let $V=\bC \cdot E_2\oplus M_{\ast}$ be the $\bC$-vector space spanned by all modular forms and the Eisenstein series $E_2$, and let $\bC\langle V\rangle$ be the shuffle algebra on $V$. Our main result is the following.
\begin{intthm}[Theorem \ref{thm:main} below]
The $QM_{\ast}$-linear morphism
\begin{align}
\varphi^{QM}: QM_{\ast} \otimes_{\bC} \bC\langle V\rangle &\rightarrow \cI^{QM}\notag\\
[f_1|\ldots|f_n] &\mapsto I(f_1,\ldots,f_n;\tau)
\end{align}
is an isomorphism of $QM_{\ast}$-algebras.
\end{intthm}
A similar result holds for the $M_{\ast}$-subalgebra $\cI^M$ of $\cI^{QM}$ of iterated integrals of modular forms (cf. Theorem \ref{thm:mainmod}).\footnote{After this paper has been submitted for publication, the author learned that, in the case of iterated integrals of modular forms, a very similar result has also been proved by Brown (cf. \cite{Brown:EquivIterIntEis}, Proposition 4.4), using a slightly different method.}
The surjectivity of $\varphi^{QM}$ can be reduced to the fact that every quasimodular form can be written uniquely as a polynomial in $n$-th derivatives of modular forms and the Eisenstein series $E_2$ (cf. \cite{123}, Proposition 20). The proof of injectivity is more elaborate and amounts to showing that iterated integrals of modular forms and the Eisenstein series $E_2$ are linearly independent over $QM_{\ast}$. It extends a result of \cite{LMS} which dealt with iterated integrals of Eisenstein series. In both cases, the key is to use a general result on linear independence of iterated integrals \cite{DDMS}. It would be interesting to prove similar results for quasimodular forms for congruence subgroups.

The Milnor--Moore theorem \cite{MM} states that if $k$ has characteristic zero, then $k\langle V\rangle$ is isomorphic to a polynomial algebra (usually in infinitely many variables). Fixing a (totally ordered) basis $\cB$ of $V$, Radford \cite{Radford} has given explicit generators of $k\langle V\rangle$ in terms of Lyndon words on $\cB$ (cf. Section \ref{sec:4}). Using this, we get the following theorem.
\begin{intthm}[Theorem \ref{thm:polynomial} below]
Let $\cB$ be a basis of $\bC \cdot E_2 \oplus M_{\ast}$. We have a natural isomorphism
\begin{equation} \label{eqn:canisom}
\cI^{QM} \cong QM_{\ast}[Lyn(\cB^*)],
\end{equation}
where the right hand side is the polynomial $QM_{\ast}$-algebra on the set $Lyn(\cB^*)$ of Lyndon words of $\cB$.
\end{intthm}
Again, a similar result holds for $\cI^M$. Since $QM_{\ast}$ has an explicit basis given by monomials in the Eisenstein series $E_2$, $E_4$ and $E_6$, the isomorphism \eqref{eqn:canisom} can be made completely explicit, and may be viewed as an analog of the isomorphism $QM_{\ast} \cong \bC[E_2,E_4,E_6]$ \cite{KanekoZagier}.

Finally, we note that classically, integrals of modular forms play an important role in Eichler--Shimura theory, where they give rise to group-cocycles (say for $\SL_2(\bZ)$ or more generally for some congruence subgroup thereof) with values in homogeneous polynomials. This has been generalized by Manin \cite{Man}, and later by Brown \cite{Brown:MMV} and Hain \cite{Hain:HodgeDeRham}, who attach certain non-abelian cocycles to iterated integrals of modular forms. Although it is not the main focus of this article, in the appendix we show how one can attach cocycles to quasimodular forms (for $\SL_2(\bZ)$), partly since we found no mention of this in the literature. On the other hand, we leave the definition and study of cocycles attached to iterated integrals of quasimodular forms for future investigation.

The plan of the paper is as follows. In Section \ref{sec:2}, we collect the necessary background on quasimodular forms and their iterated integrals. In Section \ref{sec:3}, we prove a linear independence result for iterated integrals of quasimodular forms. This result is then put to use in Section \ref{sec:4}, where the main results are proved.
In the appendix, we discuss the above-mentioned generalization of the classical Eichler--Shimura theory to quasimodular forms for $\SL_2(\bZ)$.

{\bf Acknowledgments:} Very many thanks to Pierre Lochak for bringing the article \cite{DDMS} to the author's attention. Also, many thanks to Francis Brown, Erik Panzer and the referees for corrections as well as very helpful suggestions and to Don Zagier for inspiring discussions on the appendix. The results of this paper were found while the author was a PhD student at Universit\"at Hamburg under the supervision of Ulf K\"uhn.
\section{Preliminaries} \label{sec:2}

Throughout the paper, all modular and quasimodular forms will be for $\SL_2(\bZ)$.
We fix some notation. Let $\fH=\{z \in \bC \, \vert \, \im(z)>0 \}$ be the upper half-plane with canonical coordinate $\tau$.
For every $k \in \bZ$, we have a group action of $\SL_2(\bZ)$ on the set of all functions $f: \fH \rightarrow \bC$ (not necessarily holomorphic), defined by $(\gamma,f) \mapsto f\vert_k\gamma$, where
\begin{equation}
(f\vert_k\gamma)(\tau):=(c\tau+d)^{-k}f\left(\frac{a\tau+b}{c\tau+d} \right).
\end{equation}
For fixed $\tau \in \fH$, we also define a map $X: \SL_2(\bZ) \rightarrow \bC$ by $X(\gamma)=\frac{1}{2\pi i}\frac{c}{c\tau+d}$. Note that $X$ has infinite, and thus Zariski dense, image.

\subsection{Recap of modular forms}
Denote by $M_k$ the space of modular forms of weight $k \in \bZ$. By definition, these are the holomorphic functions $f: \fH \rightarrow \bC$, which satisfy $f\vert_k\gamma=f$ for all $\gamma \in \SL_2(\bZ)$, and which are ``holomorphic at the cusp''. The latter condition means that in the Fourier expansion $f(\tau)=\sum_{n \in \bZ}a_nq^n$ (which exists since for $\gamma=\left(\begin{smallmatrix}1&1\\0&1\end{smallmatrix} \right) \in \SL_2(\bZ)$, the condition $f\vert_k\gamma=f$ is just $f(\tau+1)=f(\tau)$ for all $\tau$), all $a_n=0$ for $n<0$. Examples of modular forms include the Eisenstein series
\begin{equation}
E_{2k}(\tau)=1-\frac{4k}{B_{2k}}\sum_{n=1}^{\infty}n^{2k-1}\frac{q^n}{1-q^n}=1-\frac{4k}{B_{2k}}\sum_{n=1}^{\infty}\left( \sum_{d\vert n}d^{2k-1} \right)q^n,
\end{equation}
which is a modular form of weight $2k$, for $k \geq 2$ (the $B_{2k}$ are Bernoulli numbers). The $\bC$-vector space of all modular forms $M_{\ast}$ is a graded (for the weight) $\bC$-algebra $M_{\ast}=\bigoplus_{k \in \bZ}M_k$, which is well-known to be isomorphic to the polynomial algebra $\bC[E_4,E_6]$. Proofs of all these facts and much more on modular forms can be found for example in \cite{123}.
\subsection{Quasimodular forms}
Quasimodular forms are a generalization of modular forms, which have first been introduced in \cite{KanekoZagier} (see also \cite{BlochOk}, \S 3 and \cite{123}, \S 5.3). The definition we give here is due to W. Nahm\footnote{Cf. \cite{123}, Section 5.3.} and is also used for example in \cite{MartinRoyer}.
\begin{dfn}\label{dfn:qmod}
Let $k,p \in \bZ$ with $p \geq 0$. A \textit{quasimodular form} of weight $k$ and depth $\leq p$ is a function $f:\fH \rightarrow \bC$ with the following property: there exist holomorphic functions $f_r: \fH \rightarrow \bC$, for $0 \leq r \leq p$, which have Fourier expansions $\sum_{n=0}^{\infty}a_nq^n$, such that
\begin{equation} \label{eqn:qmod}
(f\vert_k\gamma)(\tau)=\sum_{r=0}^pf_r(\tau)X(\gamma)^r, \quad \mbox{for all }\gamma \in \SL_2(\bZ).
\end{equation}
We denote by $QM^{\leq p}_k$ the $\bC$-vector space of quasimodular forms of weight $k$ and depth $\leq p$ and set 
\begin{equation}
QM_k:=\bigcup_{p\geq 0}QM^{\leq p}_k, \quad QM_{\ast}:=\bigoplus_{k\in\bZ}QM_k.
\end{equation}
\end{dfn}
\begin{rmk} \label{rmk:unique}
\begin{enumerate}
\item[\rm (i)]
It is clear from the definition that, if $f_1 \in QM_{k_1}^{\leq p_1}$, $f_2 \in QM_{k_2}^{\leq p_2}$, then $f_1f_2 \in QM_{k_1+k_2}^{\leq p_1+p_2}$. In other words, $QM_{\ast}$ is a graded (for the weight) and filtered (for the depth) $\bC$-algebra.
\item[\rm (ii)]
Using that $X$ is Zariski dense, it is easy to see that the functions $f_r(\tau)$ are uniquely determined by $f(\tau)$. Also, applying \eqref{eqn:qmod} with $\gamma=\left(\begin{smallmatrix}1&0\\0&1\end{smallmatrix}\right)$, we see that $f_0(\tau)=f(\tau)$. In particular, every quasimodular form is holomorphic on $\fH$ and at the cusp.
\end{enumerate}
\end{rmk}

Every modular form is a quasimodular form of depth zero, more precisely, $M_k=QM_k^{\leq 0}$. An example of a quasimodular form, which is not modular is the Eisenstein series of weight two
$
E_2(\tau)=1-24\sum_{n=1}^{\infty}n\frac{q^n}{1-q^n},
$
which transforms as
\begin{equation} \label{eqn:E2trans}
(E_2\vert_2\gamma)(\tau)=E_2(\tau)+12X(\gamma)=E_2(\tau)-\frac{6i}{\pi}\frac{c}{c\tau+d},
\end{equation}
for all $\gamma \in \SL_2(\bZ)$. In particular, $E_2 \in QM_2^{\leq 1} \setminus M_2$.

The following proposition recalls basic properties of $QM_{\ast}$ that will be of use later.
\begin{prop} \label{prop:qmod}
\begin{enumerate}
\item[\rm (i)]
The $\bC$-algebra $QM_{\ast}$ is closed under the differential operator $D:=\frac{1}{2\pi i}\frac{d}{d\tau}=q\frac{d}{dq}$. More precisely, for $f$ quasimodular of weight $k$ and depth $\leq p$, we have
\begin{equation} \label{eqn:diff}
(D(f)\vert_{k+2}\gamma)(\tau)=\sum_{r=0}^{p+1}(D(f_r)(\tau)+(k-r+1)f_{r-1}(\tau))X(\gamma)^r.
\end{equation}
In particular, $D(QM_k^{\leq p}) \subset QM_{k+2}^{\leq p+1}$ for all $k,p \in \bZ$.
\bigskip

\item[\rm (ii)]
We have
\begin{equation}
QM_k=\begin{cases}\{0\}, & \mbox{if $k<0$} \\ \bC \cdot E_2, & \mbox{if $k=2$}\\ D(QM_{k-2}) \oplus M_k &\mbox{else.}\end{cases}
\end{equation}
In particular, $QM_{\ast}=\bC \cdot E_2 \oplus D(QM_{\ast}) \oplus M_{\ast}$, and
\begin{equation}
QM_{\ast}\cong\bC[E_2,E_4,E_6]
\end{equation}
as graded $\bC$-algebras.
\end{enumerate}
\end{prop}
\begin{prf}
For (i), simply apply $D$ to both sides of \eqref{eqn:qmod}. The first equality in (ii) follows from \cite{123}, Proposition 20.(iii), and the isomorphism $QM_{\ast} \cong \bC[E_2,E_4,E_6]$ is essentially a consequence of this, but can also be proved independently (cf. \cite{BlochOk}, Proposition 3.5.(ii)).
\end{prf}
\begin{rmk} \label{rmk:qmodfun}
Relaxing the condition in the definition of quasimodular forms that every $f_r$ be a holomorphic function, one can define the notion of \textit{weakly quasimodular form} of weight $k$ and depth $\leq p$ as a meromorphic function $f: \fH \rightarrow \bC$ satisfying \eqref{eqn:qmod}, but where the functions $f_r(\tau)$ are only required to be meromorphic on $\fH$ and have Fourier series of the form $\sum_{n=-M}^{\infty}a_nq^n$ ($f_r$ is ``meromorphic at the cusp''). As in the case of quasimodular forms, one shows easily that the functions $f_r(\tau)$ are uniquely determined by $f(\tau)$ (cf. Remark \ref{rmk:unique}). Moreover, Proposition \ref{prop:qmod}.(i) generalizes straightforwardly to weakly quasimodular forms.
\end{rmk}
We end this subsection with a short lemma, for which we couldn't find a suitable reference. Denote by $\Delta=\frac{1}{1728}(E_4^3-E_6^2)$ Ramanujan's cusp form of weight $12$.
\begin{lem} \label{lem:E2diff}
Let $g \in QM_{\ast} \setminus \{0\}$ and $\alpha \in \bC$ such that
\begin{equation} \label{eqn:E2diff}
D(g)=(\alpha E_2) \cdot g.
\end{equation}
Then $\alpha$ is a non-negative integer, and $g=\beta \Delta^{\alpha}$ for some $\beta \in \bC \setminus \{0\}$.
\end{lem}
\begin{prf}
Let $g=\sum_{n=0}^{\infty}a_nq^n$, so that $D(g)=\sum_{n=0}^{\infty}na_nq^n$. Comparing coefficients on both sides of \eqref{eqn:E2diff} yields that $\alpha$ equals the smallest integer $m \geq 0$ such that $a_m \neq 0$. On the other hand,
$
\frac{D(\Delta)}{\Delta}=E_2
$
(cf. \cite{123}, proof of Proposition 7), and from the chain rule $\frac{D(\Delta^{\alpha})}{\Delta^{\alpha}}=\alpha E_2$, which gives the result.
\end{prf}

\subsection{Iterated integrals on the upper half-plane} \label{ssec:2.2}
Iterated integrals of modular forms have been considered first by Manin (for cusp forms) \cite{Man}, and later by Brown (in general) \cite{Brown:MMV}. They are generalizations of the classical Eichler integrals \cite{Eichler,Lang}
\begin{equation} \label{eqn:EichlerInt}
\int_{\tau}^{i\infty}f(z)z^m\dd z, \quad m=0,\ldots,k-2
\end{equation}
where $f$ is a cusp form of weight $k$. Extending \eqref{eqn:EichlerInt} to a general modular form poses the problem of logarithmic divergences, which arise from the constant term in the Fourier series of $f$. A procedure for regularizing such integrals is described in \cite{Brown:MMV}, and we borrow it to define iterated integrals of quasimodular forms. Since it is perhaps not so well-known, we give some details, for the convenience of the reader.

Let $W \subset \cO(\fH)$ be the $\bC$-subalgebra of holomorphic functions  $f: \fH \rightarrow \bC$, which have an everywhere convergent Fourier series $f(\tau)=\sum_{n=0}^{\infty}a_nq^n$ with $q=e^{2\pi i\tau}$. Note that $QM_{\ast}\subset W$. For $f(\tau) \in W$, let $f^{\infty}=a_0$, and $f^0(\tau)=f(\tau)-f^{\infty}=\sum_{n=1}^{\infty}a_nq^n$.
Let $\bC\langle W\rangle$ (sometimes denoted by $T^c(W)$) be the shuffle algebra \cite{Reu}, i.e. the graded dual of the tensor algebra $T(W)=\bigoplus_{k \geq 0} W^{\otimes n}$ on $W$, where the grading is by the length of tensors. Elements of $(W^{\otimes n})^{\vee}$ will be written using bar notation $[f_1|f_2|\ldots|f_n]$, and a general element of $\bC\langle W\rangle$ is a $\bC$-linear combination of those. The product on $\bC\langle W\rangle$ is the shuffle product $\shuffle$, which is defined on the basic elements by
\begin{equation} \label{eqn:shuffleproduct}
[f_1|\ldots|f_r] \shuffle [f_{r+1}|\ldots|f_{r+s}]=\sum_{\sigma \in \Sigma_{r,s}}[f_{\sigma(1)}|\ldots|f_{\sigma(r+s)}],
\end{equation}
where $\Sigma_{r,s}$ denotes the set of all the permutations on the set $\{1,\ldots,r+s\}$ such that $\sigma^{-1}(1)<\ldots<\sigma^{-1}(r)$ and $\sigma^{-1}(r+1)<\ldots<\sigma^{-1}(r+s)$.

Define a $\bC$-linear map $R: \bC\langle W\rangle \rightarrow \bC\langle W\rangle$ by the formula
\begin{equation}
R[f_1|\ldots|f_n]=\sum_{i=0}^n(-1)^{n-i}[f_1|\ldots|f_i] \shuffle [f_n^{\infty}|\ldots|f_{i+1}^{\infty}].
\end{equation}
Following \cite{Brown:MMV}, Section 4, we make the following definition.
\begin{dfn} \label{dfn:maindef}
For $f_1,\ldots,f_n \in W$, define their regularized iterated integral
\begin{equation} \label{eqn:maindef}
I(f_1,\ldots,f_n;\tau):=(2\pi i)^n\sum_{i=0}^n(-1)^{n-i}\int_{\tau}^{i\infty}R[f_1|\ldots|f_i]\int_0^{\tau}[f_n^{\infty}|\ldots|f_{i+1}^{\infty}],
\end{equation}
where $\displaystyle \int\limits_a^b[f_1|\ldots|f_n]:=\int\limits_{0 \leq t_1 \leq \ldots\leq t_n\leq 1}(\gamma_a^b)^*(f_1(\tau_1)\dd\tau_1)\ldots (\gamma_a^b)^*(f_n(\tau_n)\dd\tau_n)$ denotes the usual iterated integral along the straight line path $\gamma_a^b$ from $a$ to $b$.
\end{dfn}
\begin{rmk}
Using the change of variables $\tau \mapsto q=e^{2\pi i\tau}$, it is easy to see that $I(f_1,\ldots,f_n;\tau) \in W[\log(q)]$, where $\log(q):=2\pi i\tau$. By the same token, if all of the $f_i$ have rational Fourier coefficients, then $I(f_1,\ldots,f_n;\tau)$ will also have rational coefficients, as a series in $q$ and $\log(q)$.
\end{rmk}
\begin{prop}
The functions $I(f_1,\ldots,f_n;\tau)$ satisfy the following properties.
\bigskip

\begin{enumerate}
\item[\rm (i)]
The product of any two of them is given by the shuffle product:
\begin{equation} \label{eqn:shuffle}
I(f_1,\ldots,f_r;\tau)I(f_{r+1},\ldots,f_{r+s};\tau)=\sum_{\sigma \in \Sigma_{r,s}}I(f_{\sigma(1)},\ldots,f_{\sigma(r+s)};\tau).
\end{equation}
\smallskip

\item[\rm (ii)]
They satisfy the differential equation
\begin{equation} \label{eqn:diffeq}
\frac{1}{2\pi i}\frac{d}{d\tau}\Bigr\vert_{\tau=\tau_0}I(f_1,\ldots,f_n;\tau)=-f_1(\tau_0)I(f_2,\ldots,f_n;\tau_0).
\end{equation}
\smallskip

\item[\rm(iii)]
We have the integration by parts formulas
\begin{align}
I(f_1,\ldots,f_i,D(g),f_{i+1},\ldots,f_n;\tau)&=I(f_1,\ldots,f_i,gf_{i+1},\ldots,f_n;\tau)\notag\\
&-I(f_1,\ldots,f_ig,f_{i+1},\ldots,f_n;\tau), \label{eqn:intparts}
\end{align}
as well as
\begin{equation}
I(D(g),f_2,\ldots,f_n;\tau)=I(gf_2,f_3,\ldots,f_n;\tau)-g(\tau)I(f_2,\ldots,f_n;\tau),
\end{equation}
and
\begin{equation}
I(f_1,\ldots,f_{n-1},D(g);\tau)=g(i\infty)I(f_1,\ldots,f_{n-1};\tau)-I(f_1,\ldots,f_{n-1}g;\tau).
\end{equation}
\end{enumerate}
\end{prop}
\begin{prf}
Using the definition \eqref{eqn:maindef}, all of these follow from the analogous properties for usual iterated integrals (cf. e.g. \cite{Hai}). 
\end{prf}

\subsection{A criterion for linear independence of iterated integrals}
Let $\Frac(W)$ be the field of fractions of the $\bC$-algebra $W$ introduced in the last subsection. By the quotient rule, it is easy to see that $\Frac(W)$ is closed under $D=\frac{1}{2\pi i}\frac{d}{d\tau}$.

The following theorem is a special case of the main result of \cite{DDMS}.
\begin{thm} \label{thm:DDMS}
Let $\cF=(f_i)_{i \in I}$ be a family of elements of $W$, and let $\cC \subset \Frac(W)$ be a subfield, which is closed under $D$ and contains $\cF$. The following are equivalent:
\begin{enumerate}
\item[\rm (i)]
The family of iterated integrals $(I(f_1,\ldots,f_n;\tau) \, \vert \, f_i \in I, \, n\geq 0)$ is linearly independent over $\cC$.
\item[\rm (ii)]
The family $\cF$ is linearly independent over $\bC$, and we have
\begin{equation} \label{eqn:trivint}
D (\cC) \cap \spn_{\bC}(\cF)=\{0\}.
\end{equation}
\end{enumerate}
\end{thm}
\begin{prf}
This is the special case of Theorem 2.1 in \cite{DDMS}, with the notation of \textit{loc.cit.}, $k=\bC$,  $(\cA,\dd)=(\Frac(\cO(\fH)),D)$, $X=\{A_{f_i} \, \vert \, f_i \in \cF \}$, $M=-\sum_{i \in I}f_iA_{f_i}$ and $S=\sum_{n \geq 0}\sum_{f_{i_1},\ldots,f_{i_n} \in S}I(f_1,\ldots,f_n;\tau)\cdot A_{f_1}\ldots A_{f_n}$. Note that it follows from \eqref{eqn:diffeq} that
\begin{equation}
D (S)=M\cdot S,
\end{equation}
as required in Theorem 2.1 of \cite{DDMS}.
\end{prf}
\begin{rmk}
Variants of Theorem \ref{thm:DDMS} have been known before (cf. \cite{Brown:thesis}, Lemma 3.6).
\end{rmk}
\section{Linear independence of iterated integrals of quasimodular forms} \label{sec:3}

In this section, we apply Theorem \ref{thm:DDMS} to deduce linear independence of a large family of iterated integrals of quasimodular forms. More precisely, our main result is the following theorem.
\begin{thm} \label{thm:linind}
Let $\cB$ be a $\bC$-linearly independent family of elements of $\bC \cdot E_2 \oplus M_{\ast}$. Then the family of iterated integrals
\begin{equation}
(I(f_1,\ldots,f_n;\tau) \, \vert \, f_i \in \cB)
\end{equation}
is linearly independent over $\Frac(QM_{\ast}) \cong \bC(E_2,E_4,E_6)$.
\end{thm}
\subsection{Two auxiliary lemmas}
For the proof of Theorem \ref{thm:linind}, we need two lemmas.
\begin{lem} \label{lem:1}
Let $f,g \in \bC[E_2,E_4,E_6]$ such that $g \neq 0$ and such that $f$ and $g$ are coprime. Assume that $D\left(\frac fg\right) \in \bC[E_2,E_4,E_6]$. Then $g=\beta \Delta^{\alpha}$ for some $\alpha \in \bZ_{\geq 0}$ and some $\beta \in \bC \setminus \{0\}$, where $\Delta:=\frac{1}{1728}(E_4^3-E_6^2)$ is Ramanujan's cusp form of weight 12.
\end{lem}
\begin{prf}
By the quotient rule, we have
\begin{equation}
D\left(\frac fg\right)=\frac{D(f)g-fD(g)}{g^2}=\frac{D(f)-f\frac{D(g)}{g}}{g}.
\end{equation}
The left hand side is contained in $\bC[E_2,E_4,E_6]$ by assumption, and since also $D(f)$ and $g$ are in $\bC[E_2,E_4,E_6]$, we have $f\frac{D(g)}{g} \in \bC[E_2,E_4,E_6]$. But then, as $f$ and $g$ have no common factor, $g$ must divide $D(g)$, i.e. there exists $h \in \bC[E_2,E_4,E_6]$ such that
\begin{equation} \label{eqn:product}
D(g)=gh.
\end{equation}
Since the operator $D: QM_{\ast} \rightarrow QM_{\ast}$ is homogeneous of weight $2$ (cf. Proposition \ref{prop:qmod}.(i)), we have $h \in QM_2$, i.e. $h=\alpha E_2$ with $\alpha \in \bC$. In other words, $g$ solves the differential equation $D(g)=(\alpha E_2)\cdot g$. But by Lemma \ref{lem:E2diff}, $\alpha$ must be a non-negative integer and $g=\beta\Delta^{\alpha}$ for some $\beta\in \bC\setminus \{0\}$.
\end{prf}
\begin{lem} \label{lem:2}
Let $f$ be a weakly quasimodular form, such that its derivative $D(f)$ is a quasimodular form. Then $f$ is a quasimodular form.
\end{lem}
\begin{prf}
It is no loss of generality to assume that $f$ is of weight $k \in \bZ$ and depth $\leq p$, where $p \geq 0$.
By the definition of weakly quasimodular forms (cf. also Remark \ref{rmk:unique}), there exist uniquely determined meromorphic functions $f_r(\tau)$, for $0\leq r\leq p$, such that
\begin{equation}
(f\vert_k\gamma)(\tau)=\sum_{r=0}^pf_r(\tau)X(\gamma)^r,
\end{equation}
for all $\gamma \in \SL_2(\bZ)$. Therefore, we only need to show that every $f_r(\tau)$ is holomorphic, including at the cusp.

To this end, by Proposition \ref{prop:qmod}.(i), we know that
\begin{equation} \label{eqn:auxalmhol}
(D(f)\vert_{k+2}\gamma)(\tau)=\sum_{r=0}^{p+1}(D(f_r)(\tau)+(k-r+1)f_{r-1}(\tau))X(\gamma)^r,
\end{equation}
and since $D(f)$ is a quasimodular form by assumption, every coefficient of \eqref{eqn:auxalmhol} is holomorphic, including at the cusp.

The constant term, with respect to $X(\gamma)$, in \eqref{eqn:auxalmhol} equals $D(f_0)(\tau)$, which is holomorphic by assumption. But a meromorphic function whose derivative is holomorphic everywhere is itself holomorphic everywhere. An easy induction argument, using that the coefficients of \eqref{eqn:auxalmhol} are holomorphic, now shows that in fact every $f_r(\tau)$ is holomorphic.
\end{prf}
\subsection{Proof of Theorem \ref{thm:linind}}
We will use the criterion of Theorem \ref{thm:DDMS} in the case where $\cC=\Frac(QM_{\ast})$ and $\cF=\cB$. Since $\cB$ is linearly independent over $\bC$ by assumption, it is enough to prove that if $h \in \Frac(QM_{\ast})$ then 
\begin{equation}
D(h)=\sum_{f \in \cB}\alpha_f f, \, \alpha_f \in \bC \quad \Rightarrow \quad \alpha_f=0, \, \mbox{ for all } f \in \cB.
\end{equation} Also, since $\cB$ spans a subspace of $\bC \cdot E_2 \oplus M_{\ast}$, it clearly suffices to prove that $D(h) \in \bC \cdot E_2 \oplus M_{\ast}$ implies that $D(h)=0$, or equivalently that $h$ is constant. Thus, the following proposition completes the proof of Theorem \ref{thm:linind}.
\begin{prop} \label{prop:technical}
Let $h \in \Frac(QM_{\ast}) \cong \bC(E_2,E_4,E_6)$, such that $D(h) \in \bC \cdot E_2\oplus M_{\ast}$. Then $h$ is constant.
\end{prop}
\begin{prf}
Write $h=\frac fg$ with $f,g \in \bC[E_2,E_4,E_6]$, $g\neq 0$ and such that $f$ and $g$ are coprime. Writing $f$ as a $\bC$-linear combination of its homogeneous components, it is enough to show the proposition for $f$ homogeneous of weight $k_f$.

First, we know from Lemma \ref{lem:1} that $g=\beta\Delta^{\alpha}$ for some $\alpha \in \bZ_{\geq 0}$ and $\beta \in \bC \setminus \{0\}$, where $\Delta$ is Ramanujan's cusp form of weight $12$. In particular, $g$ is a cusp form of weight $k_g=12\alpha$.

Since $f$ is quasimodular of weight $k_f$ and depth $\leq p$, there exist holomorphic (including at the cusp) functions $f_r(\tau)$, for $0\leq r\leq p$, such that
\begin{equation}
(f\vert_{k_f}\gamma)(\tau)=\sum_{r=0}^pf_r(\tau)X(\gamma)^r,
\end{equation}
for all $\gamma \in \SL_2(\bZ)$. Setting $h_r(\tau):=\frac{f_r}{g}(\tau)$, we also have, for $k:=k_f-k_g$
\begin{equation} \label{eqn:5}
(h\vert_k\gamma)(\tau)=\sum_{r=0}^ph_r(\tau)X(\gamma)^r.
\end{equation}
Moreover, the functions $h_r(\tau)$ are meromorphic, thus, $h$ is a weakly quasimodular form (of weight $k$ and depth $\leq p$). By assumption, $D(h)$ is a quasimodular form (necessarily of weight $k+2$ and depth $\leq p+1$), and using Lemma \ref{lem:2}, this implies that $h \in QM^{\leq p}_k$, therefore every $h_r(\tau)$ is holomorphic, including at the cusp.

Summarizing, we have seen that $h \in \Frac(QM_{\ast})$ such that $D(h) \in QM_{\ast}$ implies that $h \in QM_{\ast}$. But we even have $D(h) \in \bC \cdot E_2\oplus M_{\ast}$ by assumption, and therefore Proposition \ref{prop:qmod}.(ii) now implies that $h$ is constant, as was to be shown.
\end{prf}

\section{Iterated integrals of quasimodular forms and shuffle algebras}\label{sec:4}

We describe the $QM_{\ast}$-algebra of iterated integrals of quasimodular forms, which is the smallest algebra, which contains $QM_{\ast}$ and is closed under integration. Using the results of the last section, we show that it is canonically isomorphic to an explicit shuffle algebra. A similar result holds for the $M_{\ast}$-subalgebra of iterated integrals of modular forms.

\subsection{The algebra of iterated integrals of quasimodular forms}
\begin{dfn}
Define $\cI^{QM}$ to be the $QM_{\ast}$-module generated by all iterated integrals of quasimodular forms:
\begin{equation}
\cI^{QM}=\spn_{QM_{\ast}}\{ I(f_1,\ldots,f_n;\tau) \, \vert \, f_i \in QM_{\ast} \}.
\end{equation}
We also denote by $\cI^{QM}_n$ the $QM_{\ast}$-linear submodule, which is spanned by all of the $I(f_1,\ldots,f_r;\tau)$ with $r \leq n$.
\end{dfn}
The subspaces $\cI^{QM}_n$ define an ascending filtration $\cI^{QM}_{\bullet}$ on $\cI^{QM}$, called the length filtration (in analogy with the length filtration on iterated integrals \cite{Hai}). It follows from \eqref{eqn:shuffle} that $\cI^{QM}$ is a filtered $QM_{\ast}$-algebra. However, the length is not a grading, as shown by the next result.
\begin{prop} \label{prop:length}
Let $f_1,\ldots,f_n$ be quasimodular forms. Then
\begin{equation}
I(f_1,\ldots,f_{i-1},D(f_i),f_{i+1},\ldots,f_n;\tau) \in \cI^{QM}_{n-1}.
\end{equation}
\end{prop}
\begin{prf}
This is an immediate consequence of the integration by parts formula \eqref{eqn:intparts}.
\end{prf}

\subsection{$\cI^{QM}$ as a shuffle algebra}
We let $V$ be the $\bC$-vector space $\bC \cdot E_2 \oplus M_{\ast}$, and denote by $\bC\langle V\rangle$ the shuffle algebra on $V$ (cf. Section \ref{ssec:2.2}). Recall that this is the graded dual of the tensor algebra $T(V)$, whose grading is given by the length of tensors. Elements of $\bC\langle V\rangle$ are $\bC$-linear combination of the basic elements $[f_1|\ldots|f_n]$, and the product on $\bC\langle V\rangle$ is the shuffle product \eqref{eqn:shuffleproduct}.

The following theorem is the main result of this paper.
\begin{thm} \label{thm:main}
The $QM_{\ast}$-linear map
\begin{align} \label{eqn:varphi}
\varphi^{QM}: QM_{\ast} \otimes_{\bC} \bC\langle V\rangle &\rightarrow \cI^{QM}\\
[f_1|\ldots|f_n] &\mapsto I(f_1,\ldots,f_n;\tau) \notag
\end{align}
is an isomorphism of $QM_{\ast}$-algebras.
\end{thm}
\begin{prf}
Let $\cB$ be a basis of $V$, so that the family $([f_1|\ldots|f_n] \, \vert \, f_i \in \cB)$ is a basis of $\bC\langle V\rangle$. The injectivity of $\varphi^{QM}$ follows from the $\Frac(QM_{\ast})$-linear independence of the family
\begin{equation} \label{eqn:family}
\cF=(I(f_1,\ldots,f_n;\tau) \, \vert \, f_i \in \cB),
\end{equation}
which is a consequence of Theorem \ref{thm:linind}.

In order to obtain the surjectivity, we need to prove that the family \eqref{eqn:family} generates $\cI^{QM}$. To this end, we prove inductively that for every $n\geq 0$, we have $\cI^{QM}_n \subset \spn_{QM_{\ast}}\cF$. 
The case $n=0$ is trivial. Now let $n\geq 1$ and assume that for every $r \leq n-1$, we have $\cI^{QM}_r \subset \spn_{QM_{\ast}}\cF$. Given quasimodular forms $f_1,\ldots,f_n$, we can write $f_i=g_i+D(h_i)$, where $g_i \in\bC \cdot E_2 \oplus M_{\ast}$ and $h_i \in D(QM_{\ast})$ by Proposition \ref{prop:qmod}.(ii). Then by linearity
\begin{align}
I(f_1,\ldots,f_n;\tau)&=I(g_1,\ldots,g_n;\tau)\notag\\&+\sum_{i=1}^nI(g_1,\ldots,g_{i-1},D(h_i),g_{i+1},\ldots,g_n)+\ldots, \label{eqn:sum}
\end{align}
where the $\ldots$ above signifies iterated integrals, which have at least two $D(h_i)$ as integrands. The first term on the right is contained in $\spn_{QM_{\ast}}\cF$, since $g_i \in \bC \cdot E_2 \oplus M_{\ast}$ for every $i$ and $\cB$ is a basis. On the other hand, all other terms in the sum \eqref{eqn:sum} are iterated integrals, which contain at least one $D(h_i)$. By Proposition \ref{prop:length}, it thus follows that $I(f_1,\ldots,f_n;\tau) \equiv I(g_1,\ldots,g_n;\tau) \mod \cI^{QM}_{n-1}$, and we conclude using the induction hypothesis. Finally, it is clear that $\varphi^{QM}$ is a homomorphism of algebras, since both sides of \eqref{eqn:varphi} are endowed with the shuffle product.
\end{prf}
\subsection{The algebra of iterated integrals of modular forms}
In this section, we study the subalgebra $\cI^M$ of $\cI^{QM}$, generated by iterated integrals of modular forms.
\begin{dfn} \label{dfn:algitermod}
Define $\cI^M$ to be the $M_{\ast}$-module generated by all iterated integrals of modular forms:
\begin{equation}
\cI^{M}=\spn_{M_{\ast}}\{ I(f_1,\ldots,f_n;\tau) \, \vert \, f_i \in M_{\ast} \}.
\end{equation}
\end{dfn}
As in the case of $\cI^{QM}$, the length of iterated integrals defines the length filtration $\cI^{M}_{\bullet}$ on $\cI^{M}$, and $\cI^M$ is a filtered $M_{\ast}$-subalgebra of $\cI^{QM}$. We let $\bC\langle M_{\ast}\rangle$ be the shuffle algebra on the $\bC$-vector space $M_{\ast}$.
\begin{thm} \label{thm:mainmod}
The $M_{\ast}$-linear map
\begin{align} \label{eqn:varphi2}
\varphi^M: M_{\ast} \otimes_{\bC} \bC\langle M_{\ast}\rangle &\rightarrow \cI^{M}\\
[f_1|\ldots|f_n] &\mapsto I(f_1,\ldots,f_n;\tau) \notag
\end{align}
is an isomorphism of $M_{\ast}$-algebras.
\end{thm}
\begin{prf}
The morphism $\varphi^M$ is surjective by definition. It is also injective, since for a basis $\cB_M$ of $M_{\ast}$, the iterated integrals $I(f_1,\ldots,f_n;\tau)$ with $f_i \in \cB_M$ are linearly independent over $M_{\ast}$ by Theorem \ref{thm:linind}, as $M_{\ast} \subset \Frac(QM_{\ast})$.
\end{prf}

\subsection{A polynomial basis for $\cI^{QM}$}
Recall from Proposition \ref{prop:qmod}.(ii) that $QM_{\ast}$ is isomorphic to the polynomial algebra $\bC[E_2,E_4,E_6]$. A similar, but slightly more involved statement holds for the $QM_{\ast}$-algebra $\cI^{QM}$ of iterated integrals of quasimodular forms. Namely, $\cI^{QM}$ is a polynomial algebra over $QM_{\ast}$ in infinitely many variables, which are given by certain Lyndon words.

In the following, if $(S,<)$ is a totally ordered set, we will endow the free monoid $S^*$ on $S$ with the lexicographical order induced by $<$. Also, the \textit{length} of $w$ is simply the number of letters of $w$.
\begin{dfn}
A \textit{Lyndon word} on $S^*$ is a non-trivial word, $w \in S^* \setminus \{1\}$, such that for all factorizations $w=uv$ with $u,v \neq 1$, we have $w<v$. We denote by $Lyn(S^*)$ the set of all Lyndon words on $S^*$.
\end{dfn}
\begin{exmp}
Let $S=\{a,b\}$ with total order $a<b$. Then the Lyndon words on $S^*$ of length at most four are
\begin{equation}
a,b,ab,aab,abb,aaab,aabb,abbb.
\end{equation}
\end{exmp}
Now for a field $k$ and any set $S$, define $k\langle S\rangle$ to be the shuffle algebra on the free $k$-vector space generated by $S$. If $k$ is of characteristic zero, then by the Milnor--Moore theorem \cite{MM}, $k\langle S\rangle$ is isomorphic to a polynomial algebra (in possibly infinitely many variables). The following refinement is due to Radford.
\begin{thm}[\cite{Radford}] \label{thm:rad}
If $k$ has characteristic zero, then $k\langle S\rangle$ is freely generated, as a $k$-algebra, by the set of Lyndon words $Lyn(S^*)$. Equivalently, $k\langle S\rangle \cong k[Lyn(S^*)]$, the polynomial algebra on $Lyn(S^*)$.
\end{thm}
Returning to quasimodular forms, consider again the $\bC$-vector space $V=\bC \cdot E_2 \oplus M_{\ast}$, and let $\cB=\cup_{k\geq 0}\cB_k$ be the homogeneous basis of $V$, given by $\cB_k=\{E_4^aE_6^b \, \vert \, 4a+6b=k \}$ for $k \neq 2$, and $\cB_2=\{E_2\}$. The basis $\cB$ can be ordered for the lexicographical order as follows: if $E_4^aE_6^b,E_4^{a'}E_6^{b'} \in \cB_k$, then
\begin{equation}
E_4^aE_6^b < E_4^{a'}E_6^{b'} :\Leftrightarrow a<a', \mbox{ or } a=a', \mbox{ and } b< b',
\end{equation}
and if $f \in \cB_k$, $g \in \cB_{k'}$ with $k<k'$, then $f<g$.

Now, since for $f_1,\ldots,f_n \in \cB$, the iterated integrals $I(f_1,\ldots,f_n;\tau)$ are linearly independent over $QM_{\ast}$ (by Theorem \ref{thm:linind}), we can canonically identify the set of all $I(f_1,\ldots,f_n;\tau)$ with the free monoid $\cB^*$, and order $\cB^*$ for the lexicographical ordering induced from the order on $\cB$ above. The next result is a formal consequence of Theorems \ref{thm:main}, \ref{thm:mainmod} and \ref{thm:rad}. 
\begin{thm} \label{thm:polynomial}
The elements of $Lyn(\cB^*)$ are algebraically independent over $QM_{\ast}$ and we have a natural isomorphism of $QM_{\ast}$-algebras
\begin{equation}
QM_{\ast}[Lyn(\cB^*)] \cong \cI^{QM},
\end{equation}
which is filtered for the length, where the left hand side is the polynomial $QM_{\ast}$-algebra on $Lyn(\cB^*)$. Explicitly, the isomorphism maps an element $w=f_1\ldots f_n \in Lyn(\cB^*)$ to the iterated integral $I(f_1,\ldots,f_n;\tau)$. Similarly, we have a natural isomorphism of $M_{\ast}$-algebras
\begin{equation}
M_{\ast}[Lyn(\cB^*_M)] \cong \cI^M,
\end{equation}
where $\cB_M=\cB \setminus \{E_2\}$.
\end{thm}
\begin{exmp}
The following table gives all elements of $Lyn(\cB^*)$ involving iterated integrals of length at most two of quasimodular forms of total weight at most $12$. For ease of notation, we have dropped the $\tau$ from $I(f_1,\ldots,f_n;\tau)$.
\begin{center}
	{\footnotesize
		\begin{tabular}{|c|c|c|c|c|c|c|c|}
			\hline
			\diagbox{Weight}{Length} & 0&1&2 \\ \hline
			$0$ & --- & $I(1)$ & --- \\ \hline
			$2$ & --- & $I(E_2)$ & --- \\ \hline
			$4$ & --- & $I(E_4)$ & $I(1,E_4)$\\ \hline
			$6$ & --- & $I(E_6)$ & $I(1,E_6)$, $I(E_2,E_4)$\\ \hline
			$8$ & --- & $I(E_4^2)$ & $I(1,E_4^2)$, $I(E_2,E_6)$\\ \hline
			$10$& --- & $I(E_4E_6)$ & $I(1,E_4E_6)$, $I(E_2,E_4^2)$, $I(E_4,E_6)$\\ \hline
			$12$& --- & $I(E_4^3)$, $I(E_6^2)$ & $I(1,E_4^3)$, $I(1,E_6^2)$, $I(E_2,E_4E_6)$, $I(E_4,E_4^2)$\\ \hline 
		\end{tabular}
	}
\end{center}
\bigskip

Also, the list of all elements of $Lyn(\cB^*)$ consisting of iterated integrals of length at most three of quasimodular forms of total weight $12$ is given by
\begin{align}
\{&I(E_4^3),\, I(E_6^2),\,  I(1,E_4^3),\, I(1,E_6^2),\, I(E_2,E_4E_6),\, I(E_4,E_4^2), \notag\\
&I(1,1,E_4^3),\, I(1,1,E_6^2),\, I(1,E_2,E_4E_6),\, I(1,E_4,E_4^2),\, I(1,E_6,E_6),\notag\\
&I(1,E_4^2,E_4),\, I(1,E_4E_6,E_2),\, I(E_2,E_2,E_4^2),\, I(E_2,E_4,E_6),\, I(E_2,E_6,E_4) \}.
\end{align}
\end{exmp}

\appendix

\section{Eichler--Shimura for quasimodular forms}

In this appendix, we show how one can attach one-cocycles to quasimodular forms. This extends the classical Eichler--Shimura theory of the cocycles attached to modular forms, and is probably well-known to the experts, but the author does not know of a suitable reference for the precise statements.

Throughout this appendix, we will freely use some elementary concepts from the cohomology of groups, for which we refer to \cite{Weibel:Hom}, Ch. 6.

\subsection{Cocycles attached to modular forms}
We begin by briefly recalling how modular forms give rise to cocycles for $\SL_2(\bZ)$. A standard reference is \cite{Lang}, Ch. VI.

For $d \geq 0$, let $\bQ[X,Y]_d$ be the $\bQ$-vector space of homogeneous polynomials in $X$ and $Y$ of degree $d$. It is a right $\SL_2(\bZ)$-module by defining
\begin{equation}
P(X,Y)\vert_\gamma=P(aX+bY,cX+dY), \quad \mbox{for }\gamma=\begin{pmatrix}a&b\\c&d\end{pmatrix} \in \SL_2(\bZ), \quad P \in \bQ[X,Y]_d.
\end{equation}
With this action, given a modular form $f$ of weight $k \geq 2$, it is straightforward to verify that the holomorphic differential one-form 
\begin{equation}
\underline{f}(\tau):=(2\pi i)^{k-1}f(\tau)(X-\tau Y)^{k-2} \dd\tau \in \Omega^1(\fH) \otimes_{\bQ} \bQ[X,Y]_{k-2}
\end{equation}
is $\SL_2(\bZ)$-invariant, where $\SL_2(\bZ)$ acts on $\fH$ in the usual way via fractional linear transformations. Fixing a base point $\tau_0$ of $\fH$ (possibly $i\infty$), it follows from the $\SL_2(\bZ)$-invariance that the function
\begin{equation}
r_{f,\tau_0}: \SL_2(\bZ) \rightarrow \bC[X,Y]_{k-2}, \quad \gamma \mapsto \int_{\tau}^{\tau_0}\underline{f}(\tau)-\left(\int_{\gamma.\tau}^{\tau_0}\underline{f}(\tau)\right)\bigg\vert_\gamma,
\end{equation}
(regularized as in Section \ref{ssec:2.2}, if $\tau_0=i\infty$) is a one-cocycle, i.e. it satisfies $r_{f,\tau_0}(\gamma_1 \gamma_2)=r_{f,\tau_0}(\gamma_1)\vert_{\gamma_2}+r_{f,\tau_0}(\gamma_2)$ for all $\gamma_1,\gamma_2 \in \SL_2(\bZ)$. Its cohomology class does not depend on $\tau_0$, and we denote this class simply by $[r_f]$.

The same construction can also be applied to the complex conjugate $\overline{\underline{f}(\tau)}:=(-2\pi i)^{k-1}\overline{f(\tau)}(X-\overline{\tau}Y)^{k-2}\dd\overline{\tau}$ of the one-form $\underline{f}(\tau)$, and we denote by $[r_{\overline{f}}]$ the resulting cohomology class.
\begin{thm}[Eichler--Shimura] \label{thm:ES}
For every $k \geq 2$, the morphism
\begin{align}
M_k \oplus \overline{S}_k &\rightarrow H^1(\SL_2(\bZ),\bQ[X,Y]_{k-2}) \otimes_{\bQ}\bC,\\
(f,\overline{g}) &\mapsto [r_f]+[r_{\overline{g}}],
\end{align}
is an isomorphism of $\bC$-vector spaces. Here, $\overline{S}_k$ denotes the complex conjugate of the $\bC$-vector space of cusp forms of weight $k$.
\end{thm}
\subsection{Cocycles for the braid group}

The fact that $r_f$ is a cocycle hinges on the modularity of $f$. In order to incorporate quasimodular forms into the picture, we need to consider instead of $\SL_2(\bZ)$ the braid group $B_3=\langle \sigma_1,\sigma_2 : \sigma_1\sigma_2\sigma_1=\sigma_2\sigma_1\sigma_2 \rangle$ on three strands. It is a central extension
\begin{equation} \label{eqn:ses}
1 \longrightarrow \bZ \longrightarrow B_3 \longrightarrow \SL_2(\bZ) \longrightarrow 1,
\end{equation}
and also the fundamental group of the quotient of $\bC^{\times} \times \fH$ by the $\SL_2(\bZ)$-action
\begin{equation}
\gamma.(z,\tau)=((c\tau+d)z,\gamma.\tau), \quad \mbox{for }\gamma=\begin{pmatrix}a&b\\c&d\end{pmatrix} \in \SL_2(\bZ),
\end{equation}
where $\SL_2(\bZ)$ acts on $\fH$ as before. We refer to \cite{Hain:Modulispaces}, \S 8, for more details and further equivalent descriptions of $B_3$.

Next, we compute the cohomology groups $H^1(B_3,\bQ[X,Y]_d)$, where $B_3$ acts on $\bQ[X,Y]_d$ via the projection $B_3\rightarrow \SL_2(\bZ)$.
\begin{prop} \label{prop:ESM}
We have canonical isomorphisms
\begin{align}
H^1(B_3,\bQ[X,Y]_d) &\cong \begin{cases}H^1(\SL_2(\bZ),\bQ[X,Y]_d), & \mbox{for }d \geq 1,\\
\bQ, & \mbox{for }d=0.
\end{cases}
\end{align}
\end{prop}
\begin{prf}
The Hochschild--Serre spectral sequence (\cite{Weibel:Hom}, Ch. 6.8.3) associated to the extension \eqref{eqn:ses} yields an exact sequence
\begin{equation} \label{eqn:ses2}
0 \rightarrow H^1(\SL_2(\bZ),\bQ[X,Y]_d) \rightarrow H^1(B_3,\bQ[X,Y]_d) \rightarrow H^1(\bZ,\bQ[X,Y]_d)^{\SL_2(\bZ)} \rightarrow 0,
\end{equation}
where we have used that $H^2(\SL_2(\bZ),\bQ[X,Y]_d)=\{0\}$, as $\SL_2(\bZ)$ has virtual cohomological dimension equal to one. The proposition now follows easily from this.
\end{prf}

\subsection{Quasimodular forms and braid group cocycles}

In light of Theorem \ref{thm:ES}, Proposition \ref{prop:ESM} suggests to attach a one-cocycle $B_3 \rightarrow \bC$ to the Eisenstein series $E_2$. Indeed, this can be done as follows.

First, the modular transformation property of $E_2$ \eqref{eqn:E2trans} implies that the differential one-form
\begin{equation} \label{eqn:diffform}
2\pi iE_2(\tau)\dd \tau-12\frac{\dd z}{z} \in \Omega^1(\bC^{\times} \times \fH)
\end{equation}
is $\SL_2(\bZ)$-invariant, i.e. it descends to the quotient $\SL_2(\bZ) \backslash (\bC^{\times} \times \fH)$. Denote by
\begin{equation}
\underline{E_2}(\xi,\tau):=\varphi^*\left(2\pi iE_2(\tau)\dd \tau-12\frac{\dd z}{z}\right)=2\pi iE_2(\tau)\dd\tau-12\dd \xi \in \Omega^1(\bC \times \fH),
\end{equation}
the pull-back of \eqref{eqn:diffform} along the universal covering map $\varphi: \bC\times \fH \rightarrow \SL_2(\bZ)\backslash (\bC^{\times}\times \fH)$. Clearly, $\underline{E_2}(\xi,\tau)$ is $B_3$-invariant and it follows that for any base point $(\xi_0,\tau_0)$ (for example $(\xi_0,\tau_0)=(0,i\infty)$), the function
\begin{align}
r_{E_2,(\xi_0,\tau_0)}: B_3 &\rightarrow \bC \notag\\
\gamma &\mapsto \int_{(\xi,\tau)}^{(\xi_0,\tau_0)}\underline{E_2}(\xi,\tau)-\left(\int_{\gamma.(\xi,\tau)}^{(\xi_0,\tau_0)}\underline{E_2}(\xi,\tau)\right)\bigg\vert_\gamma \label{eqn:E2cocycle}
\end{align}
is a well-defined cocycle (again, regularization is needed if $\tau_0=i\infty$).
\begin{rmk}
The integral $I(E_2;\tau)$ introduced in Section \ref{ssec:2.2} is actually equal to $\int_{\tau}^{i\infty} \underline{E_2}(\xi,\tau)$, where we embed $\fH$ into $\bC \times \fH$ by $\tau \mapsto (0,\tau)$. However, that embedding is not $B_3$-equivariant, and indeed the integral $I(E_2;\tau)$ does not give rise to a cocycle for $B_3$; for this, one really needs to lift the form $2\pi iE_2(\tau)\dd\tau$ to the form $\underline{E_2}(\xi,\tau)$.
\end{rmk}
Now since the cocycle $r_{E_2,(\xi_0,\tau_0)}$ is non-zero, its cohomology class (which is again independent of the choice of base point $(\xi_0,\tau_0)$) is non-trivial. The Eichler--Shimura theorem (Theorem \ref{thm:ES}) together with Proposition \ref{prop:ESM} then implies the next result.
\begin{cor} \label{cor:ESM}
For every $k \geq 2$, the morphism
\begin{align}
V_k \oplus \overline{S}_k &\rightarrow H^1(B_3,\bQ[X,Y]_{k-2}) \otimes_{\bQ}\bC,\\
(f,\overline{g}) &\mapsto [r_f]+[r_{\overline{g}}],
\end{align}
where $V:=M_{\ast} \oplus \bC\cdot E_2$, is an isomorphism of $\bC$-vector spaces.
\end{cor}
One can also attach a cocycle $r_{f,\tau_0}$ to a general quasimodular form $f \in QM_k$ of weight $k$ as follows. By Proposition \ref{prop:qmod}.(ii), we know that $f$ can be written uniquely as a $\bC$-linear combination of derivatives of modular forms and of derivatives of $E_2$. Thus, we can write
\begin{equation}
f=\sum \lambda_g \cdot D^{p_g}(g), \quad \lambda_g \in \bC, \, p_g \geq 0,
\end{equation}
where $g$ is either a modular form of weight $k-2p_g$ or $g=E_2$. Therefore, we may define $r_{f,\tau_0}: B_3 \rightarrow \bC[X,Y]_{\leq k-2}:=\bigoplus_{0\leq d\leq k-2}\bC[X,Y]_d$ by
\begin{equation}
r_{f,\tau_0}:=\sum \lambda_g \cdot r_{g,\tau_0}.
\end{equation}
Using this definition, one sees in particular that the cocycles of quasimodular forms can be expressed in terms of the cocycles attached to modular forms and to $E_2$. This is of course in line with Corollary \ref{cor:ESM}.

\begin{rmk}
In \cite{Brown:MMV,Hain:HodgeDeRham,Man}, certain non-abelian $\SL_2(\bZ)$-cocycles given in terms of iterated integrals of modular forms are studied. It would be natural to try and extend this theory to non-abelian $B_3$-cocycles attached to iterated integrals of quasimodular forms (perhaps along the lines suggested in \cite{Hain:HodgeDeRham}, \S 14), but this is beyond the scope of the present paper.
\end{rmk}

\begin{bibtex}[\jobname]

@book {Abe,
    AUTHOR = {Abe, Eiichi},
     TITLE = {Hopf algebras},
    SERIES = {Cambridge Tracts in Mathematics},
    VOLUME = {74},
      NOTE = {Translated from the Japanese by Hisae Kinoshita and Hiroko
              Tanaka},
 PUBLISHER = {Cambridge University Press, Cambridge-New York},
      YEAR = {1980},
     PAGES = {xii+284},
      ISBN = {0-521-22240-0},
   MRCLASS = {16A24 (14L17)},
  MRNUMBER = {594432},
}

@article {ABW:sunrise,
	AUTHOR = {Adams, Luise and Bogner, Christian and Weinzierl, Stefan},
	TITLE = {The two-loop sunrise graph with arbitrary masses},
	JOURNAL = {J. Math. Phys.},
	FJOURNAL = {Journal of Mathematical Physics},
	VOLUME = {54},
	YEAR = {2013},
	NUMBER = {5},
	PAGES = {052303, 18},
	ISSN = {0022-2488},
	MRCLASS = {81Q30 (81S40)},
	MRNUMBER = {3098926},
	DOI = {10.1063/1.4804996},
	URL = {http://dx.doi.org/10.1063/1.4804996},
}

@article {And,
	AUTHOR = {Anderson, Greg W.},
	TITLE = {The hyperadelic gamma function},
	JOURNAL = {Invent. Math.},
	FJOURNAL = {Inventiones Mathematicae},
	VOLUME = {95},
	YEAR = {1989},
	NUMBER = {1},
	PAGES = {63--131},
	ISSN = {0020-9910},
	CODEN = {INVMBH},
	MRCLASS = {11S80 (11G20)},
	MRNUMBER = {969414},
	MRREVIEWER = {Gerd Faltings},
	DOI = {10.1007/BF01394145},
	URL = {http://dx.doi.org/10.1007/BF01394145},
}

@book {An,
    AUTHOR = {Andr{\'e}, Yves},
     TITLE = {Une introduction aux motifs (motifs purs, motifs mixtes,
              p\'eriodes)},
    SERIES = {Panoramas et Synth\`eses [Panoramas and Syntheses]},
    VOLUME = {17},
 PUBLISHER = {Soci\'et\'e Math\'ematique de France, Paris},
      YEAR = {2004},
     PAGES = {xii+261},
      ISBN = {2-85629-164-3},
   MRCLASS = {14F42 (11J91 14C25 19E15)},
  MRNUMBER = {2115000 (2005k:14041)},
MRREVIEWER = {Luca Barbieri Viale},
}

@article {Ap,
	AUTHOR = {Ap{\'e}ry, Roger},
	 TITLE = {Irrationalit{\'e} de {$\zeta(2)$} et {$\zeta(3)$}},
   JOURNAL = {Ast{\'e}risque},
    VOLUME = {61},
      YEAR = {1979},
     PAGES = {11--13},
}
	 
@book {B,
    AUTHOR = {Bailey, W. N.},
     TITLE = {Generalized hypergeometric series},
    SERIES = {Cambridge Tracts in Mathematics and Mathematical Physics, No.
              32},
 PUBLISHER = {Stechert-Hafner, Inc., New York},
      YEAR = {1964},
     PAGES = {v+108},
   MRCLASS = {33.20 (40.00)},
  MRNUMBER = {0185155 (32 \#2625)},
}

@article {BR,
    AUTHOR = {Ball, Keith and Rivoal, Tanguy},
     TITLE = {Irrationalit\'e d'une infinit\'e de valeurs de la fonction
              z\^eta aux entiers impairs},
   JOURNAL = {Invent. Math.},
  FJOURNAL = {Inventiones Mathematicae},
    VOLUME = {146},
      YEAR = {2001},
    NUMBER = {1},
     PAGES = {193--207},
      ISSN = {0020-9910},
   MRCLASS = {11J72 (11M06)},
  MRNUMBER = {1859021},
MRREVIEWER = {F. Beukers},
       DOI = {10.1007/s002220100168},
       URL = {http://dx.doi.org/10.1007/s002220100168},
}

@article {BannaiKings,
	AUTHOR = {Bannai, Kenichi and Kings, Guido},
	TITLE = {{$p$}-adic elliptic polylogarithm, {$p$}-adic {E}isenstein
		series and {K}atz measure},
	JOURNAL = {Amer. J. Math.},
	FJOURNAL = {American Journal of Mathematics},
	VOLUME = {132},
	YEAR = {2010},
	NUMBER = {6},
	PAGES = {1609--1654},
	ISSN = {0002-9327},
	MRCLASS = {11G55 (11F85 11S80 14F30 14F42)},
	MRNUMBER = {2766179},
	MRREVIEWER = {Jan Nekov{\'a}{\v{r}}},
}

@article {BKT,
    AUTHOR = {Bannai, Kenichi and Kobayashi, Shinichi and Tsuji, Takeshi},
     TITLE = {On the de {R}ham and {$p$}-adic realizations of the elliptic
              polylogarithm for {CM} elliptic curves},
   JOURNAL = {Ann. Sci. \'Ec. Norm. Sup\'er. (4)},
  FJOURNAL = {Annales Scientifiques de l'\'Ecole Normale Sup\'erieure.
              Quatri\`eme S\'erie},
    VOLUME = {43},
      YEAR = {2010},
    NUMBER = {2},
     PAGES = {185--234},
      ISSN = {0012-9593},
   MRCLASS = {11G55 (11G15 14F30 14G10)},
  MRNUMBER = {2662664 (2011g:11125)},
MRREVIEWER = {Jan Nekov{\'a}{\v{r}}},
}

@unpublished{BS,
   author = {{Baumard}, S. and {Schneps}, L.},
   title = "{On the derivation representation of the fundamental {L}ie algebra of mixed elliptic motives}",
   note = {arXiv:1510.05549},
   year = 2015,
}

@incollection {BeiLev,
    AUTHOR = {Be{\u\i}linson, A. and Levin, A.},
     TITLE = {The elliptic polylogarithm},
 BOOKTITLE = {Motives ({S}eattle, {WA}, 1991)},
    SERIES = {Proc. Sympos. Pure Math.},
    VOLUME = {55},
     PAGES = {123--190},
 PUBLISHER = {Amer. Math. Soc., Providence, RI},
      YEAR = {1994},
   MRCLASS = {11G05 (11G09 11G40 14H52 19F27)},
  MRNUMBER = {1265553},
MRREVIEWER = {J. Browkin},
}
		
@article {BB,
	AUTHOR = {Belkale, Prakash and Brosnan, Patrick},
	TITLE = {Matroids, motives, and a conjecture of {K}ontsevich},
	JOURNAL = {Duke Math. J.},
	FJOURNAL = {Duke Mathematical Journal},
	VOLUME = {116},
	YEAR = {2003},
	NUMBER = {1},
	PAGES = {147--188},
	ISSN = {0012-7094},
	CODEN = {DUMJAO},
	MRCLASS = {14G10 (05A15 05B35 05C05 11G25 14M12 81T18)},
	MRNUMBER = {1950482},
	MRREVIEWER = {Timothy Y. Chow},
	DOI = {10.1215/S0012-7094-03-11615-4},
	URL = {http://dx.doi.org/10.1215/S0012-7094-03-11615-4},
}

@book {Bloch,
    AUTHOR = {Bloch, Spencer},
     TITLE = {Higher regulators, algebraic {$K$}-theory, and zeta functions
              of elliptic curves},
    SERIES = {CRM Monograph Series},
    VOLUME = {11},
 PUBLISHER = {American Mathematical Society, Providence, RI},
      YEAR = {2000},
     PAGES = {x+97},
      ISBN = {0-8218-2114-8},
   MRCLASS = {11G55 (11G40 11R70 14G10 19F27)},
  MRNUMBER = {1760901 (2001i:11082)},
MRREVIEWER = {Jan Nekov{\'a}{\v{r}}},
}

@article {BlochOk,
	AUTHOR = {Bloch, Spencer and Okounkov, Andrei},
	TITLE = {The character of the infinite wedge representation},
	JOURNAL = {Adv. Math.},
	FJOURNAL = {Advances in Mathematics},
	VOLUME = {149},
	YEAR = {2000},
	NUMBER = {1},
	PAGES = {1--60},
	ISSN = {0001-8708},
	CODEN = {ADMTA4},
	MRCLASS = {11F22 (17B66)},
	MRNUMBER = {1742353},
	MRREVIEWER = {Mirko Primc},
	DOI = {10.1006/aima.1999.1845},
	URL = {http://dx.doi.org/10.1006/aima.1999.1845},
}

@article {BEK,
	AUTHOR = {Bloch, Spencer and Esnault, H{\'e}l{\`e}ne and Kreimer, Dirk},
	TITLE = {On motives associated to graph polynomials},
	JOURNAL = {Comm. Math. Phys.},
	FJOURNAL = {Communications in Mathematical Physics},
	VOLUME = {267},
	YEAR = {2006},
	NUMBER = {1},
	PAGES = {181--225},
	ISSN = {0010-3616},
	CODEN = {CMPHAY},
	MRCLASS = {81Q30 (11M41 14F25 81T18)},
	MRNUMBER = {2238909},
	MRREVIEWER = {Matilde Marcolli},
	DOI = {10.1007/s00220-006-0040-2},
	URL = {http://dx.doi.org/10.1007/s00220-006-0040-2},
}

@article{BlochVanhove,
	author         = "Bloch, Spencer and Vanhove, Pierre",
	title          = "{The elliptic dilogarithm for the sunset graph}",
	journal        = "J. Number Theory",
	volume         = "148",
	pages          = "328--364",
	year           = "2015",
	eprint         = "1309.5865",
	archivePrefix  = "arXiv",
	primaryClass   = "hep-th",
	SLACcitation   = "
}

@article {MZVmine,
    AUTHOR = {Bl{\"u}mlein, J. and Broadhurst, D. J. and Vermaseren, J. A.
              M.},
     TITLE = {The multiple zeta value data mine},
   JOURNAL = {Comput. Phys. Comm.},
  FJOURNAL = {Computer Physics Communications. An International Journal and
              Program Library for Computational Physics and Physical
              Chemistry},
    VOLUME = {181},
      YEAR = {2010},
    NUMBER = {3},
     PAGES = {582--625},
      ISSN = {0010-4655},
     CODEN = {CPHCBZ},
   MRCLASS = {11M32 (11Y70)},
  MRNUMBER = {2578167 (2011a:11163)},
MRREVIEWER = {Zhonghua Li},
       DOI = {10.1016/j.cpc.2009.11.007},
       URL = {http://dx.doi.org/10.1016/j.cpc.2009.11.007},
}

@book {Bou,
    AUTHOR = {Bourbaki, N.},
     TITLE = {\'{E}l\'ements de math\'ematique. {F}asc. {XXXIV}. {G}roupes
              et alg\`ebres de {L}ie. {C}hapitre {IV}: {G}roupes de
              {C}oxeter et syst\`emes de {T}its. {C}hapitre {V}: {G}roupes
              engendr\'es par des r\'eflexions. {C}hapitre {VI}: syst\`emes
              de racines},
    SERIES = {Actualit\'es Scientifiques et Industrielles, No. 1337},
 PUBLISHER = {Hermann, Paris},
      YEAR = {1968},
     PAGES = {288 pp. (loose errata)},
   MRCLASS = {22.50 (17.00)},
  MRNUMBER = {0240238 (39 \#1590)},
MRREVIEWER = {G. B. Seligman},
}

@article {BG,
	AUTHOR = {Borwein, Jonathan M. and Girgensohn, Roland},
	TITLE = {Evaluation of triple {E}uler sums},
	JOURNAL = {Electron. J. Combin.},
	FJOURNAL = {Electronic Journal of Combinatorics},
	VOLUME = {3},
	YEAR = {1996},
	NUMBER = {1},
	PAGES = {Research Paper 23, approx.\ 27 pp.\},
		ISSN = {1077-8926},
		MRCLASS = {11M99 (11Y60)},
		MRNUMBER = {1401442},
		MRREVIEWER = {Bruce C. Berndt},
		URL = {http://www.combinatorics.org/Volume_3/Abstracts/v3i1r23.html},
	}

@article {BK,
    AUTHOR = {Broadhurst, D. J. and Kreimer, D.},
     TITLE = {Association of multiple zeta values with positive knots via
              {F}eynman diagrams up to {$9$} loops},
   JOURNAL = {Phys. Lett. B},
  FJOURNAL = {Physics Letters. B},
    VOLUME = {393},
      YEAR = {1997},
    NUMBER = {3-4},
     PAGES = {403--412},
      ISSN = {0370-2693},
     CODEN = {PYLBAJ},
   MRCLASS = {11M41 (11Z05 57M25 81T18)},
  MRNUMBER = {1435933 (98g:11101)},
MRREVIEWER = {Louis H. Kauffman},
       DOI = {10.1016/S0370-2693(96)01623-1},
       URL = {http://dx.doi.org/10.1016/S0370-2693(96)01623-1},
}

@article {BMS,
	AUTHOR = {Broedel, Johannes and Matthes, Nils and Schlotterer, Oliver},
	TITLE = {Relations between elliptic multiple zeta values and a special
		derivation algebra},
	JOURNAL = {J. Phys. A},
	FJOURNAL = {Journal of Physics. A. Mathematical and Theoretical},
	VOLUME = {49},
	YEAR = {2016},
	NUMBER = {15},
	PAGES = {155203, 49},
	ISSN = {1751-8113},
	MRCLASS = {11M32 (33E05)},
	MRNUMBER = {3479125},
	DOI = {10.1088/1751-8113/49/15/155203},
	URL = {http://dx.doi.org/10.1088/1751-8113/49/15/155203},
}

@misc{mine,
	title = {Elliptic multiple zeta value datamine},
	howpublished = {\url{https://tools.aei.mpg.de/emzv/index.html}},
	note = {Accessed: 2016-09-17}
}

@article {BMMS,
	AUTHOR = {Broedel, Johannes and Mafra, Carlos R. and Matthes, Nils and
		Schlotterer, Oliver},
	TITLE = {Elliptic multiple zeta values and one-loop superstring
		amplitudes},
	JOURNAL = {J. High Energy Phys.},
	FJOURNAL = {Journal of High Energy Physics},
	YEAR = {2015},
	NUMBER = {7},
	PAGES = {112, front matter+41},
	ISSN = {1126-6708},
	MRCLASS = {83C47 (83E30)},
	MRNUMBER = {3383100},
	MRREVIEWER = {Farhang Loran},
}

@unpublished{BMRS,
	author = {Broedel, Johannes and Matthes, Nils and Richter, Gregor and Schlotterer, Oliver},
	title = {Twisted elliptic multiple zeta values and one-loop superstring amplitudes},
	year = {2017},
	note = {arXiv:1704.03449},
}

@article {BSS,
    AUTHOR = {Broedel, Johannes and Schlotterer, Oliver and Stieberger,
              Stephan},
     TITLE = {Polylogarithms, multiple zeta values and superstring
              amplitudes},
   JOURNAL = {Fortschr. Phys.},
  FJOURNAL = {Fortschritte der Physik. Progress of Physics},
    VOLUME = {61},
      YEAR = {2013},
    NUMBER = {9},
     PAGES = {812--870},
      ISSN = {0015-8208},
   MRCLASS = {81T30 (11M32 33B30)},
  MRNUMBER = {3104459},
MRREVIEWER = {Giuseppe Nardelli},
       DOI = {10.1002/prop.201300019},
       URL = {http://dx.doi.org/10.1002/prop.201300019},
}

@article{BSST,
      author         = "Broedel, Johannes and Schlotterer, Oliver and Stieberger,
                        Stephan and Terasoma, Tomohide",
      title          = "{All order $\alpha^{\prime}$-expansion of superstring
                        trees from the Drinfeld associator}",
      journal        = "Phys. Rev.",
      volume         = "D89",
      year           = "2014",
      number         = "6",
      pages          = "066014",
      doi            = "10.1103/PhysRevD.89.066014",
      eprint         = "1304.7304",
      archivePrefix  = "arXiv",
      primaryClass   = "hep-th",
      reportNumber   = "DAMTP-2013-23, AEI-2013-195, MPP-2013-120",
      SLACcitation   = "
}

@article {Brown:thesis,
	AUTHOR = {Brown, Francis},
	TITLE = {Multiple zeta values and periods of moduli spaces
		{$\overline{\mathscr M}_{0,n}$}},
	JOURNAL = {Ann. Sci. \'Ec. Norm. Sup\'er. (4)},
	FJOURNAL = {Annales Scientifiques de l'\'Ecole Normale Sup\'erieure.
		Quatri\`eme S\'erie},
	VOLUME = {42},
	YEAR = {2009},
	NUMBER = {3},
	PAGES = {371--489},
	ISSN = {0012-9593},
	MRCLASS = {32G20 (11G55 11M32 14F42 14H10 19F27 33B30)},
	MRNUMBER = {2543329},
	MRREVIEWER = {Wadim Zudilin},
}

@incollection {Colombia,
	AUTHOR = {Brown, Francis},
	TITLE = {Iterated integrals in quantum field theory},
	BOOKTITLE = {Geometric and topological methods for quantum field theory},
	PAGES = {188--240},
	PUBLISHER = {Cambridge Univ. Press, Cambridge},
	YEAR = {2013},
	MRCLASS = {81S40 (81T18 81T40)},
	MRNUMBER = {3098088},
	MRREVIEWER = {Roberto Quezada},
}

@incollection {Brown:Decomposition,
	AUTHOR = {Brown, Francis},
	TITLE = {On the decomposition of motivic multiple zeta values},
	BOOKTITLE = {Galois-{T}eichm\"uller theory and arithmetic geometry},
	SERIES = {Adv. Stud. Pure Math.},
	VOLUME = {63},
	PAGES = {31--58},
	PUBLISHER = {Math. Soc. Japan, Tokyo},
	YEAR = {2012},
	MRCLASS = {11M32 (13B05 16T15)},
	MRNUMBER = {3051238},
	MRREVIEWER = {Antanas Laurin{\v{c}}ikas},
}

@article {Brown:MTM,
    AUTHOR = {Brown, Francis},
     TITLE = {Mixed {T}ate motives over {$\mathbb Z$}},
   JOURNAL = {Ann. of Math. (2)},
  FJOURNAL = {Annals of Mathematics. Second Series},
    VOLUME = {175},
      YEAR = {2012},
    NUMBER = {2},
     PAGES = {949--976},
      ISSN = {0003-486X},
   MRCLASS = {11S20 (11M32 14F42)},
  MRNUMBER = {2993755},
MRREVIEWER = {Pierre A. Lochak},
       DOI = {10.4007/annals.2012.175.2.10},
       URL = {http://dx.doi.org/10.4007/annals.2012.175.2.10},
}

@unpublished{BrownAnatomy,
     author         = "Brown, F.",
     title          = "{Anatomy of an associator}",
     year			= "{2012}",
     note			= "{Notes}"
}

@unpublished{Brown:depth,
     author         = {Brown, Francis},
     title          = {Depth-graded motivic multiple zeta values},
     year           = {2013},
     note			= {arXiv:1301.3053},
}

@unpublished{Brown:EquivIterIntEis,
	author         = {Brown, Francis},
	title          = {A class of non-holomorphic modular forms {II}: Equivariant iterated Eisenstein integrals},
	year           = {2017},
	note			= {arXiv:1708.03354},
}

@article {Brown:svzold,
	AUTHOR = {Brown, Francis },
	TITLE = {Polylogarithmes multiples uniformes en une variable},
	JOURNAL = {C. R. Math. Acad. Sci. Paris},
	FJOURNAL = {Comptes Rendus Math\'ematique. Acad\'emie des Sciences. Paris},
	VOLUME = {338},
	YEAR = {2004},
	NUMBER = {7},
	PAGES = {527--532},
	ISSN = {1631-073X},
	MRCLASS = {19F27 (11G55 11R42 19E20 33B30)},
	MRNUMBER = {2057024},
	MRREVIEWER = {Alexey A. Panchishkin},
	DOI = {10.1016/j.crma.2004.02.001},
	URL = {http://dx.doi.org/10.1016/j.crma.2004.02.001},
}

@article {Brown:svz,
    AUTHOR = {Brown, Francis},
     TITLE = {Single-valued motivic periods and multiple zeta values},
   JOURNAL = {Forum Math. Sigma},
  FJOURNAL = {Forum of Mathematics. Sigma},
    VOLUME = {2},
      YEAR = {2014},
     PAGES = {e25, 37},
      ISSN = {2050-5094},
   MRCLASS = {11M32},
  MRNUMBER = {3271288},
MRREVIEWER = {A. Perelli},
       DOI = {10.1017/fms.2014.18},
       URL = {http://dx.doi.org/10.1017/fms.2014.18},
}

@unpublished{Brown:MMV,
	author         = {Brown, Francis},
	title          = {Multiple modular values and the relative completion of the fundamental group of $\mathcal{M}_{1,1}$},
	year           = {2016},
	note			= {arXiv:1407.5167v3},
}

@unpublished{Brown:depth3,
     author         = {Brown, Francis},
     title          = {Zeta elements in depth 3 and the fundamental {L}ie algebra of a punctured elliptic curve},
     year           = {2015},
     note			= {arXiv:1504.04737},
}

@article{BrownLetter,
    author     =     {Brown, Francis},
    title     =     {{Letter to the author}},
    year     =     {2015},
    }

@unpublished{BL,
     author         = {Brown, Francis. and Levin, Andrey},
     title          = {Multiple elliptic polylogarithms},
	 note   		= {arXiv:1110.6917},
     year           = {2011},
}

@article {BS:K3,
	AUTHOR = {Brown, Francis and Schnetz, Oliver},
	TITLE = {A {K}3 in {$\phi^4$}},
	JOURNAL = {Duke Math. J.},
	FJOURNAL = {Duke Mathematical Journal},
	VOLUME = {161},
	YEAR = {2012},
	NUMBER = {10},
	PAGES = {1817--1862},
	ISSN = {0012-7094},
	CODEN = {DUMJAO},
	MRCLASS = {81Q30 (14C35 14J28)},
	MRNUMBER = {2954618},
	MRREVIEWER = {Paolo Aluffi},
	DOI = {10.1215/00127094-1644201},
	URL = {http://dx.doi.org/10.1215/00127094-1644201},
}

@incollection {Car,
    AUTHOR = {Cartier, Pierre},
     TITLE = {A primer of {H}opf algebras},
 BOOKTITLE = {Frontiers in number theory, physics, and geometry. {II}},
     PAGES = {537--615},
 PUBLISHER = {Springer, Berlin},
      YEAR = {2007},
   MRCLASS = {16W30 (01A60 05E05)},
  MRNUMBER = {2290769 (2008b:16059)},
MRREVIEWER = {Ralf Holtkamp},
       DOI = {10.1007/978-3-540-30308-4_12},
       URL = {http://dx.doi.org/10.1007/978-3-540-30308-4_12},
}

@incollection {CEE,
             AUTHOR = {Calaque, Damien and Enriquez, Benjamin and Etingof, Pavel},
              TITLE = {Universal {KZB} equations: the elliptic case},
          BOOKTITLE = {Algebra, arithmetic, and geometry: in honor of {Y}u. {I}.
                       {M}anin. {V}ol. {I}},
             SERIES = {Progr. Math.},
             VOLUME = {269},
              PAGES = {165--266},
          PUBLISHER = {Birkh\"auser Boston, Inc., Boston, MA},
               YEAR = {2009},
            MRCLASS = {32G34 (11F55 17B37 20C08 32C38)},
           MRNUMBER = {2641173 (2011k:32018)},
         MRREVIEWER = {Gwyn Bellamy},
                DOI = {10.1007/978-0-8176-4745-2_5},
                URL = {http://dx.doi.org/10.1007/978-0-8176-4745-2\_5},
}

@unpublished{Calaque,
	author			= "Calaque, D. and Gonzalez, M.",
	title			= "{On the universal twisted elliptic KZB connection}",
	note			= "{to appear}",
}

@article {Ch,
    AUTHOR = {Chen, Kuo Tsai},
     TITLE = {Iterated path integrals},
   JOURNAL = {Bull. Amer. Math. Soc.},
  FJOURNAL = {Bulletin of the American Mathematical Society},
    VOLUME = {83},
      YEAR = {1977},
    NUMBER = {5},
     PAGES = {831--879},
      ISSN = {0002-9904},
   MRCLASS = {55D35 (58A99)},
  MRNUMBER = {0454968 (56 \#13210)},
MRREVIEWER = {Jean-Michel Lemaire},
}

@article {D'HokerGreenVanhove,
	AUTHOR = {D'Hoker, Eric and Green, Michael B. and Vanhove, Pierre},
	TITLE = {On the modular structure of the genus-one type {II}
		superstring low energy expansion},
	JOURNAL = {J. High Energy Phys.},
	FJOURNAL = {Journal of High Energy Physics},
	YEAR = {2015},
	NUMBER = {8},
	PAGES = {041, front matter+66},
	ISSN = {1126-6708},
	MRCLASS = {81T30},
	MRNUMBER = {3402124},
}

@incollection {Del,
   AUTHOR = {Deligne, P.},
    TITLE = {Le groupe fondamental de la droite projective moins trois
             points},
BOOKTITLE = {Galois groups over ${\bf Q}$ ({B}erkeley, {CA}, 1987)},
   SERIES = {Math. Sci. Res. Inst. Publ.},
   VOLUME = {16},
    PAGES = {79--297},
PUBLISHER = {Springer, New York},
     YEAR = {1989},
  MRCLASS = {14G25 (11G35 11M06 11R70 14F35 19E99 19F27)},
 MRNUMBER = {1012168 (90m:14016)},
MRREVIEWER = {James Milne},
      DOI = {10.1007/978-1-4613-9649-9_3},
      URL = {http://dx.doi.org/10.1007/978-1-4613-9649-9\_3},
}

@article {Del23468,
	AUTHOR = {Deligne, Pierre},
	TITLE = {Le groupe fondamental unipotent motivique de {$\bold G_m-\mu_N$}, pour {$N=2,3,4,6$} ou {$8$}},
	JOURNAL = {Publ. Math. Inst. Hautes \'Etudes Sci.},
	FJOURNAL = {Publications Math\'ematiques. Institut de Hautes \'Etudes
		Scientifiques},
	NUMBER = {112},
	YEAR = {2010},
	PAGES = {101--141},
	ISSN = {0073-8301},
	MRCLASS = {14F42},
	MRNUMBER = {2737978},
	MRREVIEWER = {Claudio Pedrini},
	DOI = {10.1007/s10240-010-0027-6},
	URL = {http://dx.doi.org/10.1007/s10240-010-0027-6},
}

@article {DelMZV,
    AUTHOR = {Deligne, Pierre},
     TITLE = {Multiz\^etas, d'apr\`es {F}rancis {B}rown},
      NOTE = {S{\'e}minaire Bourbaki. Vol. 2011/2012. Expos{\'e}s
              1043--1058},
   JOURNAL = {Ast\'erisque},
  FJOURNAL = {Ast\'erisque},
    NUMBER = {352},
      YEAR = {2013},
     PAGES = {Exp. No. 1048, viii, 161--185},
      ISSN = {0303-1179},
      ISBN = {978-2-85629-371-3},
   MRCLASS = {11S40 (11G09 14C15 14F35)},
  MRNUMBER = {3087346},
MRREVIEWER = {Damian R{\"o}ssler},
}

@article {DG,
    AUTHOR = {Deligne, Pierre and Goncharov, Alexander B.},
     TITLE = {Groupes fondamentaux motiviques de {T}ate mixte},
   JOURNAL = {Ann. Sci. \'Ecole Norm. Sup. (4)},
  FJOURNAL = {Annales Scientifiques de l'\'Ecole Normale Sup\'erieure.
              Quatri\`eme S\'erie},
    VOLUME = {38},
      YEAR = {2005},
    NUMBER = {1},
     PAGES = {1--56},
      ISSN = {0012-9593},
     CODEN = {ASENAH},
   MRCLASS = {11G55 (14F42 14G10 19F27)},
  MRNUMBER = {2136480 (2006b:11066)},
MRREVIEWER = {Tam{\'a}s Szamuely},
       DOI = {10.1016/j.ansens.2004.11.001},
       URL = {http://dx.doi.org/10.1016/j.ansens.2004.11.001},
}

@incollection {DDMS,
    AUTHOR = {Deneufch{\^a}tel, Matthieu and Duchamp, G{\'e}rard H. E. and
              Minh, Vincel Hoang Ngoc and Solomon, Allan I.},
     TITLE = {Independence of hyperlogarithms over function fields via
              algebraic combinatorics},
 BOOKTITLE = {Algebraic informatics},
    SERIES = {Lecture Notes in Comput. Sci.},
    VOLUME = {6742},
     PAGES = {127--139},
 PUBLISHER = {Springer, Heidelberg},
      YEAR = {2011},
   MRCLASS = {Database Expansion Item},
  MRNUMBER = {2846744},
       DOI = {10.1007/978-3-642-21493-6_8},
       URL = {http://dx.doi.org/10.1007/978-3-642-21493-6_8},
}

@preamble{
   "\def\cprime{$'$} "
}
@article {Dr,
    AUTHOR = {Drinfel{\cprime}d, V. G.},
     TITLE = {On quasitriangular quasi-{H}opf algebras and on a group that
              is closely connected with {${\rm Gal}(\overline{\bf Q}/{\bf
              Q})$}},
   JOURNAL = {Algebra i Analiz},
  FJOURNAL = {Algebra i Analiz},
    VOLUME = {2},
      YEAR = {1990},
    NUMBER = {4},
     PAGES = {149--181},
      ISSN = {0234-0852},
   MRCLASS = {16W30 (17B37)},
  MRNUMBER = {1080203 (92f:16047)},
MRREVIEWER = {Ivan Penkov},
}

@article {Duhr,
	AUTHOR = {Duhr, Claude},
	TITLE = {Hopf algebras, coproducts and symbols: an application to
		{H}iggs boson amplitudes},
	JOURNAL = {J. High Energy Phys.},
	FJOURNAL = {Journal of High Energy Physics},
	YEAR = {2012},
	NUMBER = {8},
	PAGES = {043, front matter + 45},
	ISSN = {1126-6708},
	MRCLASS = {16T05 (33B30 81U99)},
	MRNUMBER = {3006955},
	MRREVIEWER = {{\^A}ngela Mestre},
}

@article {Ecalle,
    AUTHOR = {Ecalle, Jean},
     TITLE = {A{RI}/{GARI}, la dimorphie et l'arithm\'etique des
              multiz\^etas: un premier bilan},
   JOURNAL = {J. Th\'eor. Nombres Bordeaux},
  FJOURNAL = {Journal de Th\'eorie des Nombres de Bordeaux},
    VOLUME = {15},
      YEAR = {2003},
    NUMBER = {2},
     PAGES = {411--478},
      ISSN = {1246-7405},
   MRCLASS = {11M41 (11G55 19F27 33B30)},
  MRNUMBER = {2140864},
MRREVIEWER = {Alexey A. Panchishkin},
       URL = {http://jtnb.cedram.org/item?id=JTNB_2003__15_2_411_0},
}

@article {Eichler,
	AUTHOR = {Eichler, M.},
	TITLE = {Eine {V}erallgemeinerung der {A}belschen {I}ntegrale},
	JOURNAL = {Math. Z.},
	FJOURNAL = {Mathematische Zeitschrift},
	VOLUME = {67},
	YEAR = {1957},
	PAGES = {267--298},
	ISSN = {0025-5874},
	MRCLASS = {33.0X},
	MRNUMBER = {0089928},
	MRREVIEWER = {H. Cohn},
}

@article {Enriquez:EllAss,
    AUTHOR = {Enriquez, Benjamin},
     TITLE = {Elliptic associators},
   JOURNAL = {Selecta Math. (N.S.)},
  FJOURNAL = {Selecta Mathematica. New Series},
    VOLUME = {20},
      YEAR = {2014},
    NUMBER = {2},
     PAGES = {491--584},
      ISSN = {1022-1824},
   MRCLASS = {17B35 (11M32 14H10 16S30 20F36)},
  MRNUMBER = {3177926},
       DOI = {10.1007/s00029-013-0137-3},
       URL = {http://dx.doi.org/10.1007/s00029-013-0137-3},
}

@article{Eemzv,
      author        = {Benjamin Enriquez},
      title         = {Analogues elliptiques des nombres multiz{\'e}tas},
      journal 		= {Bull. Soc. Math. France},
      volume        = {144},
      number        = {3},
      year          = {2016},
      note          = {to appear. arXiv:1301.3042}
}

@unpublished{EL,
	author = {Enriquez, Benjamin and Lochak, Pierre},
	title  = {Homology of depth-graded motivic {L}ie algebras and {K}oszulity},
	year   = {2014},
	note   = {arXiv:1407.4060},
}

@incollection {Fal,
    AUTHOR = {Faltings, Gerd},
     TITLE = {Mathematics around {K}im's new proof of {S}iegel's theorem},
 BOOKTITLE = {Diophantine geometry},
    SERIES = {CRM Series},
    VOLUME = {4},
     PAGES = {173--188},
 PUBLISHER = {Ed. Norm., Pisa},
      YEAR = {2007},
   MRCLASS = {14G99 (11G30 14F20)},
  MRNUMBER = {2349654 (2009i:14029)},
}

@article {FurStab,
    AUTHOR = {Furusho, Hidekazu},
     TITLE = {The multiple zeta value algebra and the stable derivation
              algebra},
   JOURNAL = {Publ. Res. Inst. Math. Sci.},
  FJOURNAL = {Kyoto University. Research Institute for Mathematical
              Sciences. Publications},
    VOLUME = {39},
      YEAR = {2003},
    NUMBER = {4},
     PAGES = {695--720},
      ISSN = {0034-5318},
     CODEN = {KRMPBV},
   MRCLASS = {11M41 (14G32)},
  MRNUMBER = {2025460},
       URL = {http://projecteuclid.org/euclid.prims/1145476044},
}

@article {Fur,
    AUTHOR = {Furusho, Hidekazu},
     TITLE = {Double shuffle relation for associators},
   JOURNAL = {Ann. of Math. (2)},
  FJOURNAL = {Annals of Mathematics. Second Series},
    VOLUME = {174},
      YEAR = {2011},
    NUMBER = {1},
     PAGES = {341--360},
      ISSN = {0003-486X},
     CODEN = {ANMAAH},
   MRCLASS = {14G32 (11G55 11M32 16W60)},
  MRNUMBER = {2811601 (2012i:14031)},
MRREVIEWER = {Pierre A. Lochak},
       DOI = {10.4007/annals.2011.174.1.9},
       URL = {http://dx.doi.org/10.4007/annals.2011.174.1.9},
}

@incollection {Furusho:MZV,
	AUTHOR = {Furusho, Hidekazu},
	TITLE = {Multiple zeta values and {G}rothendieck-{T}eichm\"uller
		groups},
	BOOKTITLE = {Primes and knots},
	SERIES = {Contemp. Math.},
	VOLUME = {416},
	PAGES = {49--82},
	PUBLISHER = {Amer. Math. Soc., Providence, RI},
	YEAR = {2006},
	MRCLASS = {14G32 (11M41)},
	MRNUMBER = {2276136},
	MRREVIEWER = {Alexey A. Panchishkin},
	DOI = {10.1090/conm/416/07887},
	URL = {http://dx.doi.org/10.1090/conm/416/07887},
}

@incollection {GKZ,
    AUTHOR = {Gangl, Herbert and Kaneko, Masanobu and Zagier, Don},
     TITLE = {Double zeta values and modular forms},
 BOOKTITLE = {Automorphic forms and zeta functions},
     PAGES = {71--106},
 PUBLISHER = {World Sci. Publ., Hackensack, NJ},
      YEAR = {2006},
   MRCLASS = {11M41 (11F11)},
  MRNUMBER = {2208210 (2006m:11138)},
MRREVIEWER = {Hirofumi Tsumura},
       DOI = {10.1142/9789812774415_0004},
       URL = {http://dx.doi.org/10.1142/9789812774415_0004},
}

@incollection {Gonell,
	AUTHOR = {Goncharov, Alexander},
	TITLE = {Mixed elliptic motives},
	BOOKTITLE = {Galois representations in arithmetic algebraic geometry
		({D}urham, 1996)},
	SERIES = {London Math. Soc. Lecture Note Ser.},
	VOLUME = {254},
	PAGES = {147--221},
	PUBLISHER = {Cambridge Univ. Press, Cambridge},
	YEAR = {1998},
	MRCLASS = {11G40 (11G09 11G55 14C35 14G10 19F27)},
	MRNUMBER = {1696477},
	MRREVIEWER = {Jan Nekov{\'a}{\v{r}}},
	DOI = {10.1017/CBO9780511662010.005},
	URL = {http://dx.doi.org/10.1017/CBO9780511662010.005},
}

@article {GonGalSym,
		AUTHOR = {Goncharov, A. B.},
		TITLE = {Galois symmetries of fundamental groupoids and noncommutative
			geometry},
		JOURNAL = {Duke Math. J.},
		FJOURNAL = {Duke Mathematical Journal},
		VOLUME = {128},
		YEAR = {2005},
		NUMBER = {2},
		PAGES = {209--284},
		ISSN = {0012-7094},
		CODEN = {DUMJAO},
		MRCLASS = {11G55 (11G09 14C30 16W30 19E15 20F34)},
		MRNUMBER = {2140264},
		MRREVIEWER = {Matilde Marcolli},
		DOI = {10.1215/S0012-7094-04-12822-2},
		URL = {http://dx.doi.org/10.1215/S0012-7094-04-12822-2},
}

@article {GonMod,
    AUTHOR = {Goncharov, A. B.},
     TITLE = {Multiple polylogarithms, cyclotomy and modular complexes},
   JOURNAL = {Math. Res. Lett.},
  FJOURNAL = {Mathematical Research Letters},
    VOLUME = {5},
      YEAR = {1998},
    NUMBER = {4},
     PAGES = {497--516},
      ISSN = {1073-2780},
   MRCLASS = {11G55 (11F67 11R42 19E20 19F15 19F27)},
  MRNUMBER = {1653320 (2000c:11108)},
MRREVIEWER = {Alexey A. Panchishkin},
       DOI = {10.4310/MRL.1998.v5.n4.a7},
       URL = {http://dx.doi.org/10.4310/MRL.1998.v5.n4.a7},
}

@incollection {GonMZV,
    AUTHOR = {Goncharov, Alexander B.},
     TITLE = {Multiple {$\zeta$}-values, {G}alois groups, and geometry of
              modular varieties},
 BOOKTITLE = {European {C}ongress of {M}athematics, {V}ol. {I} ({B}arcelona,
              2000)},
    SERIES = {Progr. Math.},
    VOLUME = {201},
     PAGES = {361--392},
 PUBLISHER = {Birkh\"auser, Basel},
      YEAR = {2001},
   MRCLASS = {11G55 (11G40 11M41 14G35 33B30 81Q30)},
  MRNUMBER = {1905330},
MRREVIEWER = {Jan Nekov{\'a}{\v{r}}},
}

@unpublished{Gon,
   author = {{Goncharov}, A.~B.},
    title = "{Multiple polylogarithms and mixed Tate motives}",
    note  = {arXiv:math/0103059},
     year = 2001,
}

@article {GonGal,
	AUTHOR = {Goncharov, A. B.},
	TITLE = {Galois symmetries of fundamental groupoids and noncommutative
		geometry},
	JOURNAL = {Duke Math. J.},
	FJOURNAL = {Duke Mathematical Journal},
	VOLUME = {128},
	YEAR = {2005},
	NUMBER = {2},
	PAGES = {209--284},
	ISSN = {0012-7094},
	CODEN = {DUMJAO},
	MRCLASS = {11G55 (11G09 14C30 16W30 19E15 20F34)},
	MRNUMBER = {2140264},
	MRREVIEWER = {Matilde Marcolli},
	DOI = {10.1215/S0012-7094-04-12822-2},
	URL = {http://dx.doi.org/10.1215/S0012-7094-04-12822-2},
}

@article {GM,
    AUTHOR = {Goncharov, A. B. and Manin, Yu. I.},
     TITLE = {Multiple {$\zeta$}-motives and moduli spaces {$\overline{\mathscr
              M}_{0,n}$}},
   JOURNAL = {Compos. Math.},
  FJOURNAL = {Compositio Mathematica},
    VOLUME = {140},
      YEAR = {2004},
    NUMBER = {1},
     PAGES = {1--14},
      ISSN = {0010-437X},
   MRCLASS = {11G55 (11M41 14H10)},
  MRNUMBER = {2004120 (2005c:11090)},
MRREVIEWER = {Gilberto Bini},
       DOI = {10.1112/S0010437X03000125},
       URL = {http://dx.doi.org/10.1112/S0010437X03000125},
}

@article{GreenSchwarz:1982,
      author         = "Green, Michael B. and Schwarz, John H.",
      title          = "{Supersymmetric Dual String Theory. 3. Loops and
                        Renormalization}",
      journal        = "Nucl.Phys.",
      volume         = "B198",
      pages          = "441-460",
      doi            = "10.1016/0550-3213(82)90334-0",
      year           = "1982",
      reportNumber   = "CALT-68-873",
      SLACcitation   = "
}

@incollection {Hai,
    AUTHOR = {Hain, Richard},
     TITLE = {The geometry of the mixed {H}odge structure on the fundamental
              group},
 BOOKTITLE = {Algebraic geometry, {B}owdoin, 1985 ({B}runswick, {M}aine,
              1985)},
    SERIES = {Proc. Sympos. Pure Math.},
    VOLUME = {46},
     PAGES = {247--282},
 PUBLISHER = {Amer. Math. Soc., Providence, RI},
      YEAR = {1987},
   MRCLASS = {14F40 (14C30 32C40 32G20 55P62 58A14)},
  MRNUMBER = {927984 (89g:14010)},
MRREVIEWER = {Toshitake Kohno},
}

@incollection {HaiPolylog,
    AUTHOR = {Hain, Richard M.},
     TITLE = {Classical polylogarithms},
 BOOKTITLE = {Motives ({S}eattle, {WA}, 1991)},
    SERIES = {Proc. Sympos. Pure Math.},
    VOLUME = {55},
     PAGES = {3--42},
 PUBLISHER = {Amer. Math. Soc., Providence, RI},
      YEAR = {1994},
   MRCLASS = {19F99 (11G99 11R70 19D55)},
  MRNUMBER = {1265550 (94k:19002)},
MRREVIEWER = {Philippe Blanc},
}

@misc{Hain:P1,
	author = {Hain, Richard},
    title = {Lectures on the {H}odge--de {R}ham theory of the fundamental group of $\mathbb{P}^1 - \{0,1,\infty\}$},
	howpublished = {\url{https://services.math.duke.edu/~hain/aws/lectures_provisional.pdf}. Accessed: 2016-10-18},
}
    
@unpublished{HaiKZB,
	author = {{Hain}, R.},
	title = "{Notes on the Universal Elliptic KZB Equation}",
	note = {arXiv:1309.0580},
	year = 2013,
}

@incollection {Hain:HodgeDeRham,
	AUTHOR = {Hain, Richard},
	TITLE = {The {H}odge--de {R}ham theory of modular groups},
	BOOKTITLE = {Recent advances in {H}odge theory},
	SERIES = {London Math. Soc. Lecture Note Ser.},
	VOLUME = {427},
	PAGES = {422--514},
	PUBLISHER = {Cambridge Univ. Press, Cambridge},
	YEAR = {2016},
	MRCLASS = {14D07 (14Gxx 32S35 58A14)},
	MRNUMBER = {3409885},
}

@incollection {Hain:Modulispaces,
	AUTHOR = {Hain, Richard},
	TITLE = {Lectures on moduli spaces of elliptic curves},
	BOOKTITLE = {Transformation groups and moduli spaces of curves},
	SERIES = {Adv. Lect. Math. (ALM)},
	VOLUME = {16},
	PAGES = {95--166},
	PUBLISHER = {Int. Press, Somerville, MA},
	YEAR = {2011},
	MRCLASS = {14H10 (14D23 14H52)},
	MRNUMBER = {2883686},
	MRREVIEWER = {Arvid Siqveland},
}

@unpublished{HM,
	author = {{Hain}, R. and {Matsumoto}, M.},
	title = "{Universal Mixed Elliptic Motives}",
    note = {arXiv:1512.03975},
	year = 2015,
}

@article {Hoffman:MHS,
	AUTHOR = {Hoffman, Michael E.},
	TITLE = {Multiple harmonic series},
	JOURNAL = {Pacific J. Math.},
	FJOURNAL = {Pacific Journal of Mathematics},
	VOLUME = {152},
	YEAR = {1992},
	NUMBER = {2},
	PAGES = {275--290},
	ISSN = {0030-8730},
	MRCLASS = {11M06 (40A05)},
	MRNUMBER = {1141796},
	MRREVIEWER = {Bruce C. Berndt},
	URL = {http://projecteuclid.org/euclid.pjm/1102636166},
}

@article {Hoff,
    AUTHOR = {Hoffman, Michael E.},
     TITLE = {Quasi-shuffle products},
   JOURNAL = {J. Algebraic Combin.},
  FJOURNAL = {Journal of Algebraic Combinatorics. An International Journal},
    VOLUME = {11},
      YEAR = {2000},
    NUMBER = {1},
     PAGES = {49--68},
      ISSN = {0925-9899},
     CODEN = {JAOME7},
   MRCLASS = {05E05},
  MRNUMBER = {1747062},
MRREVIEWER = {Jason Fulman},
       DOI = {10.1023/A:1008791603281},
       URL = {http://dx.doi.org/10.1023/A:1008791603281},
}

@article {HubKing,
    AUTHOR = {Huber, Annette and Kings, Guido},
     TITLE = {Degeneration of {$l$}-adic {E}isenstein classes and of the
              elliptic polylog},
   JOURNAL = {Invent. Math.},
  FJOURNAL = {Inventiones Mathematicae},
    VOLUME = {135},
      YEAR = {1999},
    NUMBER = {3},
     PAGES = {545--594},
      ISSN = {0020-9910},
     CODEN = {INVMBH},
   MRCLASS = {11G18 (11G55 11R42 11R70 19F27)},
  MRNUMBER = {1669288},
MRREVIEWER = {Jan Nekov{\'a}{\v{r}}},
       DOI = {10.1007/s002220050295},
       URL = {http://dx.doi.org/10.1007/s002220050295},
}

@incollection {Iha,
    AUTHOR = {Ihara, Yasutaka},
     TITLE = {The {G}alois representation arising from {${\bf P}^1-\{0,1,\infty\}$} and {T}ate twists of even degree},
 BOOKTITLE = {Galois groups over {${\bf Q}$} ({B}erkeley, {CA}, 1987)},
    SERIES = {Math. Sci. Res. Inst. Publ.},
    VOLUME = {16},
     PAGES = {299--313},
 PUBLISHER = {Springer, New York},
      YEAR = {1989},
   MRCLASS = {11F80 (11G20 11R23 11R58 14E22 14G25)},
  MRNUMBER = {1012169},
MRREVIEWER = {Sheldon Kamienny},
       DOI = {10.1007/978-1-4613-9649-9_4},
       URL = {http://dx.doi.org/10.1007/978-1-4613-9649-9_4},
}

@inproceedings {IhICM,
    AUTHOR = {Ihara, Yasutaka},
     TITLE = {Braids, {G}alois groups, and some arithmetic functions},
 BOOKTITLE = {Proceedings of the {I}nternational {C}ongress of
              {M}athematicians, {V}ol.\ {I}, {II} ({K}yoto, 1990)},
     PAGES = {99--120},
 PUBLISHER = {Math. Soc. Japan, Tokyo},
      YEAR = {1991},
   MRCLASS = {11G09 (11R32 14E20 16W30 20F34)},
  MRNUMBER = {1159208},
MRREVIEWER = {J. Browkin},
}

@article {IKZ,
    AUTHOR = {Ihara, Kentaro and Kaneko, Masanobu and Zagier, Don},
     TITLE = {Derivation and double shuffle relations for multiple zeta
              values},
   JOURNAL = {Compos. Math.},
  FJOURNAL = {Compositio Mathematica},
    VOLUME = {142},
      YEAR = {2006},
    NUMBER = {2},
     PAGES = {307--338},
      ISSN = {0010-437X},
   MRCLASS = {11M41},
  MRNUMBER = {2218898},
MRREVIEWER = {David Bradley},
       DOI = {10.1112/S0010437X0500182X},
       URL = {http://dx.doi.org/10.1112/S0010437X0500182X},
}

@article {IO,
    AUTHOR = {Ihara, Kentaro and Ochiai, Hiroyuki},
     TITLE = {Symmetry on linear relations for multiple zeta values},
   JOURNAL = {Nagoya Math. J.},
  FJOURNAL = {Nagoya Mathematical Journal},
    VOLUME = {189},
      YEAR = {2008},
     PAGES = {49--62},
      ISSN = {0027-7630},
     CODEN = {NGMJA2},
   MRCLASS = {11M41},
  MRNUMBER = {2396583},
MRREVIEWER = {David Bradley},
       URL = {http://projecteuclid.org/euclid.nmj/1205156910},
}

@incollection {KanekoZagier,
	AUTHOR = {Kaneko, Masanobu and Zagier, Don},
	TITLE = {A generalized {J}acobi theta function and quasimodular forms},
	BOOKTITLE = {The moduli space of curves ({T}exel {I}sland, 1994)},
	SERIES = {Progr. Math.},
	VOLUME = {129},
	PAGES = {165--172},
	PUBLISHER = {Birkh\"auser Boston, Boston, MA},
	YEAR = {1995},
	MRCLASS = {11F11 (11F03)},
	MRNUMBER = {1363056},
	MRREVIEWER = {Bruce Hunt},
	DOI = {10.1007/978-1-4612-4264-2_6},
	URL = {http://dx.doi.org/10.1007/978-1-4612-4264-2_6},
}

@book {Kas,
    AUTHOR = {Kassel, Christian},
     TITLE = {Quantum groups},
    SERIES = {Graduate Texts in Mathematics},
    VOLUME = {155},
 PUBLISHER = {Springer-Verlag, New York},
      YEAR = {1995},
     PAGES = {xii+531},
      ISBN = {0-387-94370-6},
   MRCLASS = {17B37 (16W30 18D10 20F36 57M25 81R50)},
  MRNUMBER = {1321145 (96e:17041)},
MRREVIEWER = {Yu. N. Bespalov},
       DOI = {10.1007/978-1-4612-0783-2},
       URL = {http://dx.doi.org/10.1007/978-1-4612-0783-2},
}

@article {Kim,
    AUTHOR = {Kim, Minhyong},
     TITLE = {{$p$}-adic {$L$}-functions and {S}elmer varieties associated
              to elliptic curves with complex multiplication},
   JOURNAL = {Ann. of Math. (2)},
  FJOURNAL = {Annals of Mathematics. Second Series},
    VOLUME = {172},
      YEAR = {2010},
    NUMBER = {1},
     PAGES = {751--759},
      ISSN = {0003-486X},
     CODEN = {ANMAAH},
   MRCLASS = {11G15 (11G05 11G40)},
  MRNUMBER = {2680431 (2011i:11089)},
MRREVIEWER = {Francesc C. Castell{\`a}},
       DOI = {10.4007/annals.2010.172.751},
       URL = {http://dx.doi.org/10.4007/annals.2010.172.751},
}

@article {Kings:Tamagawa,
	AUTHOR = {Kings, Guido},
	TITLE = {The {T}amagawa number conjecture for {CM} elliptic curves},
	JOURNAL = {Invent. Math.},
	FJOURNAL = {Inventiones Mathematicae},
	VOLUME = {143},
	YEAR = {2001},
	NUMBER = {3},
	PAGES = {571--627},
	ISSN = {0020-9910},
	CODEN = {INVMBH},
	MRCLASS = {11G40 (11G05 11G55 14G10)},
	MRNUMBER = {1817645},
	MRREVIEWER = {Jan Nekov{\'a}{\v{r}}},
	DOI = {10.1007/s002220000115},
	URL = {http://dx.doi.org/10.1007/s002220000115},
}

@article {KZ,
    AUTHOR = {Knizhnik, V. G. and Zamolodchikov, A. B.},
     TITLE = {Current algebra and {W}ess-{Z}umino model in two dimensions},
   JOURNAL = {Nuclear Phys. B},
  FJOURNAL = {Nuclear Physics. B},
    VOLUME = {247},
      YEAR = {1984},
    NUMBER = {1},
     PAGES = {83--103},
      ISSN = {0550-3213},
     CODEN = {NUPBBO},
   MRCLASS = {81E13 (81D15)},
  MRNUMBER = {853258},
       DOI = {10.1016/0550-3213(84)90374-2},
       URL = {http://dx.doi.org/10.1016/0550-3213(84)90374-2},
}

@incollection {KZPer,
	AUTHOR = {Kontsevich, Maxim and Zagier, Don},
	TITLE = {Periods},
	BOOKTITLE = {Mathematics unlimited---2001 and beyond},
	PAGES = {771--808},
	PUBLISHER = {Springer, Berlin},
	YEAR = {2001},
	MRCLASS = {11-02 (11F67 11G40 11G55)},
	MRNUMBER = {1852188},
	MRREVIEWER = {F. Beukers},
}

@book {Lang,
	AUTHOR = {Lang, Serge},
	TITLE = {Introduction to modular forms},
	NOTE = {Grundlehren der mathematischen Wissenschaften, No. 222},
	PUBLISHER = {Springer-Verlag, Berlin-New York},
	YEAR = {1976},
	PAGES = {ix+261},
	MRCLASS = {10DXX},
	MRNUMBER = {0429740},
	MRREVIEWER = {Neal Koblitz},
}

@article {Landen,
    AUTHOR = {Landen, John},
     TITLE = {Mathematical memoirs respecting a variety of subjects},
   JOURNAL = {Nourse},
      YEAR = {1780},
} 

@book {Lap,
	AUTHOR = {Lappo-Danilevsky, J. A.},
	TITLE = {M\'emoires sur la th\'eorie des syst\`emes des \'equations
		diff\'erentielles lin\'eaires},
	PUBLISHER = {Chelsea Publishing Co., New York, N. Y.},
	YEAR = {1953},
	PAGES = {xiv+253+208+204},
	MRCLASS = {36.0X},
	MRNUMBER = {0054111},
}

@article {LM,
    AUTHOR = {Le, Thang Tu Quoc and Murakami, Jun},
     TITLE = {Kontsevich's integral for the {K}auffman polynomial},
   JOURNAL = {Nagoya Math. J.},
  FJOURNAL = {Nagoya Mathematical Journal},
    VOLUME = {142},
      YEAR = {1996},
     PAGES = {39--65},
      ISSN = {0027-7630},
     CODEN = {NGMJA2},
   MRCLASS = {57M25 (11M99)},
  MRNUMBER = {1399467 (97d:57009)},
MRREVIEWER = {Sergei K. Lando},
       URL = {http://projecteuclid.org/euclid.nmj/1118772043},
} 
 
@article {Lev,
    AUTHOR = {Levin, Andrey},
     TITLE = {Elliptic polylogarithms: an analytic theory},
   JOURNAL = {Compositio Math.},
  FJOURNAL = {Compositio Mathematica},
    VOLUME = {106},
      YEAR = {1997},
    NUMBER = {3},
     PAGES = {267--282},
      ISSN = {0010-437X},
     CODEN = {CMPMAF},
   MRCLASS = {11F37 (11F27 11G40 11R70 19F27)},
  MRNUMBER = {1457106 (98d:11048)},
MRREVIEWER = {Alexey A. Panchishkin},
       DOI = {10.1023/A:1000193320513},
       URL = {http://dx.doi.org/10.1023/A:1000193320513},
}

@unpublished{LR,
     author         = {A. Levin and G. Racinet},
     title          = {Towards multiple elliptic polylogarithms},
     note           = {arXiv:math/0703237},
     year           = {2007},
}

@unpublished{LMS,
   author = {Lochak, Pierre and Matthes, Nils and Schneps, Leila},
    title = "{Elliptic multiple zeta values and the elliptic double shuffle relations}",
    year = {2017},
    note = {arXiv:1703.09410},
}

@incollection {Man,
    AUTHOR = {Manin, Yuri I.},
     TITLE = {Iterated integrals of modular forms and noncommutative modular
              symbols},
 BOOKTITLE = {Algebraic geometry and number theory},
    SERIES = {Progr. Math.},
    VOLUME = {253},
     PAGES = {565--597},
 PUBLISHER = {Birkh\"auser Boston, Boston, MA},
      YEAR = {2006},
   MRCLASS = {11F67 (11G55 11M41)},
  MRNUMBER = {2263200 (2008a:11062)},
MRREVIEWER = {Caterina Consani},
       DOI = {10.1007/978-0-8176-4532-8_10},
       URL = {http://dx.doi.org/10.1007/978-0-8176-4532-8\_10},
}

@incollection {MartinRoyer,
	AUTHOR = {Martin, Fran{\c{c}}ois and Royer, Emmanuel},
	TITLE = {Formes modulaires et p\'eriodes},
	BOOKTITLE = {Formes modulaires et transcendance},
	SERIES = {S\'emin. Congr.},
	VOLUME = {12},
	PAGES = {1--117},
	PUBLISHER = {Soc. Math. France, Paris},
	YEAR = {2005},
	MRCLASS = {11F67 (11F11 11F25)},
	MRNUMBER = {2186573},
	MRREVIEWER = {Dominic A. Lanphier},
}


@article {Matthes:Edzv,
	AUTHOR = {Matthes, Nils},
	TITLE = {Elliptic double zeta values},
	JOURNAL = {J. Number Theory},
	FJOURNAL = {Journal of Number Theory},
	VOLUME = {171},
	YEAR = {2017},
	PAGES = {227--251},
	ISSN = {0022-314X},
	CODEN = {JNUTA9},
	MRCLASS = {11M32 (11F50)},
	MRNUMBER = {3556684},
	DOI = {10.1016/j.jnt.2016.07.010},
	URL = {http://dx.doi.org/10.1016/j.jnt.2016.07.010},
}

@incollection{Matthes:Decomposition,
	author = {{Matthes}, N.},
	title = "{Decomposition of elliptic multiple zeta values and iterated Eisenstein integrals}",
	booktitle = {Various aspects of multiple zeta value 2016},
	editor = {Hidekazu Furusho},
	series = {RIMS Kokyuroku},
	volume = {2015},
	pages = {170--183},
	publisher = {Res. Inst. Math. Sci. (RIMS), Kyoto},
	year = {2017},
	archivePrefix = "arXiv",
	eprint = {1703.09597},
}

@unpublished{Matthes:linind,
   author = {{Matthes}, N.},
    title = "{Linear independence of indefinite iterated Eisenstein integrals}",
    note = {arXiv:1601.05743},
    year = 2016,
}

@unpublished{Matthes:Metab,
	author = {{Matthes}, N.},
	title = "{The meta-abelian elliptic KZB associator and periods of Eisenstein series}",
	note = {arXiv:1608.00740},
	year = 2016,
}

@phdthesis{Matthes:Thesis,
	author     =     {Matthes, Nils},
	title     =     {{Elliptic multiple zeta values}},
	school     =     {Universit\"at Hamburg},
	year     =     {2016},
}

@article {MM,
	AUTHOR = {Milnor, John W. and Moore, John C.},
	TITLE = {On the structure of {H}opf algebras},
	JOURNAL = {Ann. of Math. (2)},
	FJOURNAL = {Annals of Mathematics. Second Series},
	VOLUME = {81},
	YEAR = {1965},
	PAGES = {211--264},
	ISSN = {0003-486X},
	MRCLASS = {55.40 (55.34)},
	MRNUMBER = {0174052},
	MRREVIEWER = {I. M. James},
	DOI = {10.2307/1970615},
	URL = {http://dx.doi.org/10.2307/1970615},
}

@book {Mum,
    AUTHOR = {Mumford, David},
     TITLE = {Tata lectures on theta. {II}},
    SERIES = {Progress in Mathematics},
    VOLUME = {43},
      NOTE = {Jacobian theta functions and differential equations,
              With the collaboration of C. Musili, M. Nori, E. Previato, M.
              Stillman and H. Umemura},
 PUBLISHER = {Birkh\"auser Boston, Inc., Boston, MA},
      YEAR = {1984},
     PAGES = {xiv+272},
      ISBN = {0-8176-3110-0},
   MRCLASS = {14K25 (14H40 32G20)},
  MRNUMBER = {742776},
MRREVIEWER = {M. Kh. Gizatullin},
       DOI = {10.1007/978-0-8176-4578-6},
       URL = {http://dx.doi.org/10.1007/978-0-8176-4578-6},
}

@article {Naka,
    AUTHOR = {Nakamura, Hiroaki},
     TITLE = {On exterior {G}alois representations associated with open
              elliptic curves},
   JOURNAL = {J. Math. Sci. Univ. Tokyo},
  FJOURNAL = {The University of Tokyo. Journal of Mathematical Sciences},
    VOLUME = {2},
      YEAR = {1995},
    NUMBER = {1},
     PAGES = {197--231},
      ISSN = {1340-5705},
   MRCLASS = {11G05 (11F80 11G16 14H30)},
  MRNUMBER = {1348028},
MRREVIEWER = {Yasutaka Ihara},
}

@article{Oprisa,
	author         = "Oprisa, D. and Stieberger, S.",
	title          = "{Six gluon open superstring disk amplitude, multiple
		hypergeometric series and Euler-Zagier sums}",
	year           = "2005",
	eprint         = "hep-th/0509042",
	archivePrefix  = "arXiv",
	primaryClass   = "hep-th",
	reportNumber   = "LMU-ASC-07-05, MPP-2005-12",
	SLACcitation   = "
}

@article {Pellarin,
		AUTHOR = {Pellarin, Federico},
		TITLE = {La structure diff\'erentielle de l'anneau des formes
			quasi-modulaires pour {${\rm SL}_2(\bZ)$}},
		JOURNAL = {J. Th\'eor. Nombres Bordeaux},
		FJOURNAL = {Journal de Th\'eorie des Nombres de Bordeaux},
		VOLUME = {18},
		YEAR = {2006},
		NUMBER = {1},
		PAGES = {241--264},
		ISSN = {1246-7405},
		MRCLASS = {11J91 (11F11)},
		MRNUMBER = {2245884},
		MRREVIEWER = {Fran{\c{c}}ois Martin},
		URL = {http://jtnb.cedram.org/item?id=JTNB_2006__18_1_241_0},
}

@MastersThesis{Pollack:Thesis,
    author     =     {Pollack, Aaron},
    title     =     {{Relations between derivations arising from modular forms}},
    school     =     {Duke University},
    year     =     {2009},
    }

@article {Rac,
    AUTHOR = {Racinet, Georges},
     TITLE = {Doubles m\'elanges des polylogarithmes multiples aux racines
              de l'unit\'e},
   JOURNAL = {Publ. Math. Inst. Hautes \'Etudes Sci.},
  FJOURNAL = {Publications Math\'ematiques. Institut de Hautes \'Etudes
              Scientifiques},
    NUMBER = {95},
      YEAR = {2002},
     PAGES = {185--231},
      ISSN = {0073-8301},
   MRCLASS = {11G55 (11M41)},
  MRNUMBER = {1953193 (2004c:11117)},
MRREVIEWER = {Jan Nekov{\'a}{\v{r}}},
       DOI = {10.1007/s102400200004},
       URL = {http://dx.doi.org/10.1007/s102400200004},
}

@article {Radford,
	AUTHOR = {Radford, David E.},
	TITLE = {A natural ring basis for the shuffle algebra and an
		application to group schemes},
	JOURNAL = {J. Algebra},
	FJOURNAL = {Journal of Algebra},
	VOLUME = {58},
	YEAR = {1979},
	NUMBER = {2},
	PAGES = {432--454},
	ISSN = {0021-8693},
	CODEN = {JALGA4},
	MRCLASS = {16A24 (14L15)},
	MRNUMBER = {540649 (80k:16016)},
	MRREVIEWER = {E. J. Taft},
	DOI = {10.1016/0021-8693(79)90171-6},
	URL = {http://dx.doi.org/10.1016/0021-8693(79)90171-6},
}

@incollection {Rama,
    AUTHOR = {Ramakrishnan, Dinakar},
     TITLE = {Analogs of the {B}loch-{W}igner function for higher
              polylogarithms},
 BOOKTITLE = {Applications of algebraic {$K$}-theory to algebraic geometry
              and number theory, {P}art {I}, {II} ({B}oulder, {C}olo.,
              1983)},
    SERIES = {Contemp. Math.},
    VOLUME = {55},
     PAGES = {371--376},
 PUBLISHER = {Amer. Math. Soc., Providence, RI},
      YEAR = {1986},
   MRCLASS = {11R70 (18F25 33A70)},
  MRNUMBER = {862642},
MRREVIEWER = {Glenn Stevens},
       DOI = {10.1090/conm/055.1/862642},
       URL = {http://dx.doi.org/10.1090/conm/055.1/862642},
}

@article {Ree,
    AUTHOR = {Ree, Rimhak},
     TITLE = {Lie elements and an algebra associated with shuffles},
   JOURNAL = {Ann. of Math. (2)},
  FJOURNAL = {Annals of Mathematics. Second Series},
    VOLUME = {68},
      YEAR = {1958},
     PAGES = {210--220},
      ISSN = {0003-486X},
   MRCLASS = {17.00 (20.00)},
  MRNUMBER = {0100011 (20 \#6447)},
MRREVIEWER = {P. M. Cohn},
}

@book {Reu,
    AUTHOR = {Reutenauer, Christophe},
     TITLE = {Free {L}ie algebras},
    SERIES = {London Mathematical Society Monographs. New Series},
    VOLUME = {7},
      NOTE = {Oxford Science Publications},
 PUBLISHER = {The Clarendon Press, Oxford University Press, New York},
      YEAR = {1993},
     PAGES = {xviii+269},
      ISBN = {0-19-853679-8},
   MRCLASS = {17-02 (05-02 17B05)},
  MRNUMBER = {1231799 (94j:17002)},
MRREVIEWER = {Hartmut Laue},
}

@article {SS,
    AUTHOR = {Schlotterer, O. and Stieberger, S.},
     TITLE = {Motivic multiple zeta values and superstring amplitudes},
   JOURNAL = {J. Phys. A},
  FJOURNAL = {Journal of Physics. A. Mathematical and Theoretical},
    VOLUME = {46},
      YEAR = {2013},
    NUMBER = {47},
     PAGES = {475401, 37},
      ISSN = {1751-8113},
   MRCLASS = {81T30 (14E18)},
  MRNUMBER = {3126883},
MRREVIEWER = {Jihye Sofia Seo},
       DOI = {10.1088/1751-8113/46/47/475401},
       URL = {http://dx.doi.org/10.1088/1751-8113/46/47/475401},
}

@article {Schneps:Poisson,
	AUTHOR = {Schneps, Leila},
	TITLE = {On the {P}oisson bracket on the free {L}ie algebra in two
		generators},
	JOURNAL = {J. Lie Theory},
	FJOURNAL = {Journal of Lie Theory},
	VOLUME = {16},
	YEAR = {2006},
	NUMBER = {1},
	PAGES = {19--37},
	ISSN = {0949-5932},
	MRCLASS = {17B63 (11F11 17B01)},
	MRNUMBER = {2196410},
	MRREVIEWER = {Maurizio Martino},
}

@unpublished{Leila,
	author = {{Schneps}, L.},
	title = "{Elliptic multiple zeta values, Grothendieck-Teichm\"uller and mould theory}",
	year = 2015,
	note = {arXiv:1506.09050},
}

@article {Schnetz,
	AUTHOR = {Schnetz, Oliver},
	TITLE = {Graphical functions and single-valued multiple polylogarithms},
	JOURNAL = {Commun. Number Theory Phys.},
	FJOURNAL = {Communications in Number Theory and Physics},
	VOLUME = {8},
	YEAR = {2014},
	NUMBER = {4},
	PAGES = {589--675},
	ISSN = {1931-4523},
	MRCLASS = {81Q30},
	MRNUMBER = {3318386},
	DOI = {10.4310/CNTP.2014.v8.n4.a1},
	URL = {http://dx.doi.org/10.4310/CNTP.2014.v8.n4.a1},
}

@book {Ser,
    AUTHOR = {Serre, Jean-Pierre},
     TITLE = {Repr\'esentations lin\'eaires des groupes finis},
   EDITION = {revised},
 PUBLISHER = {Hermann, Paris},
      YEAR = {1978},
     PAGES = {182},
      ISBN = {2-7056-5630-8},
   MRCLASS = {20-01 (20C99)},
  MRNUMBER = {543841 (80f:20001)},
}

@book {SerLie,
    AUTHOR = {Serre, Jean-Pierre},
     TITLE = {Lie algebras and {L}ie groups},
    SERIES = {Lecture Notes in Mathematics},
    VOLUME = {1500},
      NOTE = {1964 lectures given at Harvard University,
              Corrected fifth printing of the second (1992) edition},
 PUBLISHER = {Springer-Verlag, Berlin},
      YEAR = {2006},
     PAGES = {viii+168},
      ISBN = {978-3-540-55008-2; 3-540-55008-9},
   MRCLASS = {17-01 (22-01)},
  MRNUMBER = {2179691},
}

@book {Serre:CohGal,
	AUTHOR = {Serre, Jean-Pierre},
	TITLE = {Cohomologie galoisienne},
	SERIES = {Lecture Notes in Mathematics},
	VOLUME = {5},
	EDITION = {Fifth},
	PUBLISHER = {Springer-Verlag, Berlin},
	YEAR = {1994},
	PAGES = {x+181},
	ISBN = {3-540-58002-6},
	MRCLASS = {12G05 (11R34)},
	MRNUMBER = {1324577},
	DOI = {10.1007/BFb0108758},
	URL = {http://dx.doi.org/10.1007/BFb0108758},
}

@article {Sou,
    AUTHOR = {Soud{\`e}res, Ismael},
     TITLE = {Motivic double shuffle},
   JOURNAL = {Int. J. Number Theory},
  FJOURNAL = {International Journal of Number Theory},
    VOLUME = {6},
      YEAR = {2010},
    NUMBER = {2},
     PAGES = {339--370},
      ISSN = {1793-0421},
   MRCLASS = {11G55 (11M32 14H10)},
  MRNUMBER = {2646761},
MRREVIEWER = {Matilde N. Lal{\'{\i}}n},
       DOI = {10.1142/S1793042110002995},
       URL = {http://dx.doi.org/10.1142/S1793042110002995},
}

@article {Stan,
    AUTHOR = {Stanley, Richard P.},
     TITLE = {Invariants of finite groups and their applications to
              combinatorics},
   JOURNAL = {Bull. Amer. Math. Soc. (N.S.)},
  FJOURNAL = {American Mathematical Society. Bulletin. New Series},
    VOLUME = {1},
      YEAR = {1979},
    NUMBER = {3},
     PAGES = {475--511},
      ISSN = {0273-0979},
     CODEN = {BAMOAD},
   MRCLASS = {20C15 (05A15 13H10 14H10 15A72 51F15)},
  MRNUMBER = {526968},
MRREVIEWER = {Ralph Strebel},
       DOI = {10.1090/S0273-0979-1979-14597-X},
       URL = {http://dx.doi.org/10.1090/S0273-0979-1979-14597-X},
}

@article {Stieberger,
	AUTHOR = {Stieberger, S.},
	TITLE = {Closed superstring amplitudes, single-valued multiple zeta
		values and the {D}eligne associator},
	JOURNAL = {J. Phys. A},
	FJOURNAL = {Journal of Physics. A. Mathematical and Theoretical},
	VOLUME = {47},
	YEAR = {2014},
	NUMBER = {15},
	PAGES = {155401, 15},
	ISSN = {1751-8113},
	MRCLASS = {81T30},
	MRNUMBER = {3191676},
	DOI = {10.1088/1751-8113/47/15/155401},
	URL = {http://dx.doi.org/10.1088/1751-8113/47/15/155401},
}

@article {Ter,
    AUTHOR = {Terasoma, Tomohide},
     TITLE = {Mixed {T}ate motives and multiple zeta values},
   JOURNAL = {Invent. Math.},
  FJOURNAL = {Inventiones Mathematicae},
    VOLUME = {149},
      YEAR = {2002},
    NUMBER = {2},
     PAGES = {339--369},
      ISSN = {0020-9910},
     CODEN = {INVMBH},
   MRCLASS = {11G55 (11M41 19F27)},
  MRNUMBER = {1918675},
MRREVIEWER = {Jan Nekov{\'a}{\v{r}}},
       DOI = {10.1007/s002220200218},
       URL = {http://dx.doi.org/10.1007/s002220200218},
}

@incollection {Tera,
    AUTHOR = {Terasoma, Tomohide},
     TITLE = {Geometry of multiple zeta values},
 BOOKTITLE = {International {C}ongress of {M}athematicians. {V}ol. {II}},
     PAGES = {627--635},
 PUBLISHER = {Eur. Math. Soc., Z\"urich},
      YEAR = {2006},
   MRCLASS = {14C30 (11M41 14F42)},
  MRNUMBER = {2275614 (2008e:14009)},
MRREVIEWER = {Jan Nekov{\'a}{\v{r}}},
}

@article {Tsu,
    AUTHOR = {Tsunogai, Hiroshi},
     TITLE = {On some derivations of {L}ie algebras related to {G}alois
              representations},
   JOURNAL = {Publ. Res. Inst. Math. Sci.},
  FJOURNAL = {Kyoto University. Research Institute for Mathematical
              Sciences. Publications},
    VOLUME = {31},
      YEAR = {1995},
    NUMBER = {1},
     PAGES = {113--134},
      ISSN = {0034-5318},
     CODEN = {KRMPBV},
   MRCLASS = {11G05 (11R32 14H30 17B40)},
  MRNUMBER = {1317526},
MRREVIEWER = {Douglas L. Ulmer},
       DOI = {10.2977/prims/1195164794},
       URL = {http://dx.doi.org/10.2977/prims/1195164794},
}

@unpublished{Wald,
	author               = {Waldschmidt, Michel},
	howpublished         = {Notes of Lectures given at Chennai IMSc, April 2011},
	note                 = {\url{http://www.math.jussieu.fr/~miw/articles/pdf/MZV2011IMSc.pdf}. Accessed 2016-10-18},
	title                = {Lectures on {M}ultiple {Z}eta {V}alues, IMSC 2011},
}

@book {Was,
    AUTHOR = {Wasow, Wolfgang},
     TITLE = {Asymptotic expansions for ordinary differential equations},
      NOTE = {Reprint of the 1976 edition},
 PUBLISHER = {Dover Publications, Inc., New York},
      YEAR = {1987},
     PAGES = {x+374},
      ISBN = {0-486-65456-7},
   MRCLASS = {34-02 (34E05)},
  MRNUMBER = {919406},
}

@book {Wat,
    AUTHOR = {Waterhouse, William C.},
     TITLE = {Introduction to affine group schemes},
    SERIES = {Graduate Texts in Mathematics},
    VOLUME = {66},
 PUBLISHER = {Springer-Verlag, New York-Berlin},
      YEAR = {1979},
     PAGES = {xi+164},
      ISBN = {0-387-90421-2},
   MRCLASS = {14-01 (14Lxx 20G99)},
  MRNUMBER = {547117 (82e:14003)},
MRREVIEWER = {M. Kh. Gizatullin},
}

@book {Weibel:Hom,
	AUTHOR = {Weibel, Charles A.},
	TITLE = {An introduction to homological algebra},
	SERIES = {Cambridge Studies in Advanced Mathematics},
	VOLUME = {38},
	PUBLISHER = {Cambridge University Press, Cambridge},
	YEAR = {1994},
	PAGES = {xiv+450},
	ISBN = {0-521-43500-5; 0-521-55987-1},
	MRCLASS = {18-01 (16-01 17-01 20-01 55Uxx)},
	MRNUMBER = {1269324},
	MRREVIEWER = {Kenneth A. Brown},
	DOI = {10.1017/CBO9781139644136},
	URL = {http://dx.doi.org/10.1017/CBO9781139644136},
}

@book {W,
    AUTHOR = {Weil, Andr{\'e}},
     TITLE = {Elliptic functions according to {E}isenstein and {K}ronecker},
      NOTE = {Ergebnisse der Mathematik und ihrer Grenzgebiete, Band 88},
 PUBLISHER = {Springer-Verlag, Berlin-New York},
      YEAR = {1976},
     PAGES = {ii+93},
      ISBN = {3-540-07422-8},
   MRCLASS = {10DXX (01A55 10-03)},
  MRNUMBER = {0562289 (58 \#27769a)},
MRREVIEWER = {S. Chowla},
}

@book {WhitWat,
    AUTHOR = {Whittaker, E. T. and Watson, G. N.},
     TITLE = {A course of modern analysis},
    SERIES = {Cambridge Mathematical Library},
      NOTE = {An introduction to the general theory of infinite processes
              and of analytic functions; with an account of the principal
              transcendental functions,
              Reprint of the fourth (1927) edition},
 PUBLISHER = {Cambridge University Press, Cambridge},
      YEAR = {1996},
     PAGES = {vi+608},
      ISBN = {0-521-58807-3},
   MRCLASS = {01A75 (30-01 33-01)},
  MRNUMBER = {1424469},
       DOI = {10.1017/CBO9780511608759},
       URL = {http://dx.doi.org/10.1017/CBO9780511608759},
}

@book {Wild,
    AUTHOR = {Wildeshaus, J{\"o}rg},
     TITLE = {Realizations of polylogarithms},
    SERIES = {Lecture Notes in Mathematics},
    VOLUME = {1650},
 PUBLISHER = {Springer-Verlag, Berlin},
      YEAR = {1997},
     PAGES = {xii+343},
      ISBN = {3-540-62460-0},
   MRCLASS = {11G18 (11G09 11G40 14D07 14G35 19F27)},
  MRNUMBER = {1482233},
MRREVIEWER = {Jan Nekov{\'a}{\v{r}}},
}

@article {ZagEll,
    AUTHOR = {Zagier, Don},
     TITLE = {The {B}loch-{W}igner-{R}amakrishnan polylogarithm function},
   JOURNAL = {Math. Ann.},
  FJOURNAL = {Mathematische Annalen},
    VOLUME = {286},
      YEAR = {1990},
    NUMBER = {1-3},
     PAGES = {613--624},
      ISSN = {0025-5831},
     CODEN = {MAANA},
   MRCLASS = {11R42 (11R70 19F27)},
  MRNUMBER = {1032949 (90k:11153)},
MRREVIEWER = {V. Kumar Murty},
       DOI = {10.1007/BF01453591},
       URL = {http://dx.doi.org/10.1007/BF01453591},
}

@article {Z,
    AUTHOR = {Zagier, Don},
     TITLE = {Periods of modular forms and {J}acobi theta functions},
   JOURNAL = {Invent. Math.},
  FJOURNAL = {Inventiones Mathematicae},
    VOLUME = {104},
      YEAR = {1991},
    NUMBER = {3},
     PAGES = {449--465},
      ISSN = {0020-9910},
     CODEN = {INVMBH},
   MRCLASS = {11F67 (11F27 11F55)},
  MRNUMBER = {1106744 (92e:11052)},
MRREVIEWER = {Rolf Berndt},
       DOI = {10.1007/BF01245085},
       URL = {http://dx.doi.org/10.1007/BF01245085},
}

@incollection {ZagPol,
    AUTHOR = {Zagier, Don},
     TITLE = {Polylogarithms, {D}edekind zeta functions and the algebraic
              {$K$}-theory of fields},
 BOOKTITLE = {Arithmetic algebraic geometry ({T}exel, 1989)},
    SERIES = {Progr. Math.},
    VOLUME = {89},
     PAGES = {391--430},
 PUBLISHER = {Birkh\"auser Boston, Boston, MA},
      YEAR = {1991},
   MRCLASS = {11R42 (11R70 19F15 19F27)},
  MRNUMBER = {1085270},
MRREVIEWER = {V. Kumar Murty},
}

@incollection {Zag,
    AUTHOR = {Zagier, Don},
     TITLE = {Values of zeta functions and their applications},
 BOOKTITLE = {First {E}uropean {C}ongress of {M}athematics, {V}ol.\ {II}
              ({P}aris, 1992)},
    SERIES = {Progr. Math.},
    VOLUME = {120},
     PAGES = {497--512},
 PUBLISHER = {Birkh\"auser, Basel},
      YEAR = {1994},
   MRCLASS = {11M41 (11F67 11G40 19F27)},
  MRNUMBER = {1341859 (96k:11110)},
MRREVIEWER = {Fernando Rodr{\'{\i}}guez Villegas},
}

@article {ZagDouble,
    AUTHOR = {Zagier, Don},
     TITLE = {Periods of modular forms, traces of {H}ecke operators, and
              multiple zeta values},
      NOTE = {Research into automorphic forms and $L$ functions (Japanese)
              (Kyoto, 1992)},
   JOURNAL = {S\=urikaisekikenky\=usho K\=oky\=uroku},
  FJOURNAL = {S\=urikaisekikenky\=usho K\=oky\=uroku},
    NUMBER = {843},
      YEAR = {1993},
     PAGES = {162--170},
   MRCLASS = {11F67 (11F25 11F72)},
  MRNUMBER = {1296720},
MRREVIEWER = {Pavel Guerzhoy},
}

@incollection {ZagDilog,
	AUTHOR = {Zagier, Don},
	TITLE = {The dilogarithm function},
	BOOKTITLE = {Frontiers in number theory, physics, and geometry. {II}},
	PAGES = {3--65},
	PUBLISHER = {Springer, Berlin},
	YEAR = {2007},
	MRCLASS = {33B30 (11G55)},
	MRNUMBER = {2290758},
	MRREVIEWER = {David Bradley},
	DOI = {10.1007/978-3-540-30308-4_1},
	URL = {http://dx.doi.org/10.1007/978-3-540-30308-4_1},
}

@incollection {123,
    AUTHOR = {Zagier, Don},
     TITLE = {Elliptic modular forms and their applications},
 BOOKTITLE = {The 1-2-3 of modular forms},
    SERIES = {Universitext},
     PAGES = {1--103},
 PUBLISHER = {Springer, Berlin},
      YEAR = {2008},
   MRCLASS = {11F11 (11-02 11E45 11F20 11F25 11F27 11F67)},
  MRNUMBER = {2409678},
MRREVIEWER = {Rainer Schulze-Pillot},
       DOI = {10.1007/978-3-540-74119-0_1},
       URL = {http://dx.doi.org/10.1007/978-3-540-74119-0_1},
}

@unpublished{Zerbini,
	author = {{Zerbini}, F.},
	title = "{Single-valued multiple zeta values in genus 1 superstring amplitudes}",
	note = {arXiv:1512.05689},
	year = {2015},
}

@article {Zud,
    AUTHOR = {Zudilin, V. V.},
     TITLE = {One of the numbers {$\zeta(5)$}, {$\zeta(7)$}, {$\zeta(9)$},
              {$\zeta(11)$} is irrational},
   JOURNAL = {Uspekhi Mat. Nauk},
  FJOURNAL = {Rossi\u\i skaya Akademiya Nauk. Moskovskoe Matematicheskoe
              Obshchestvo. Uspekhi Matematicheskikh Nauk},
    VOLUME = {56},
      YEAR = {2001},
    NUMBER = {4(340)},
     PAGES = {149--150},
      ISSN = {0042-1316},
   MRCLASS = {11J72 (11J91 11M06)},
  MRNUMBER = {1861452},
MRREVIEWER = {John H. Loxton},
       DOI = {10.1070/RM2001v056n04ABEH000427},
       URL = {http://dx.doi.org/10.1070/RM2001v056n04ABEH000427},
}

\end{bibtex}
\bibliographystyle{abbrv}
\bibliography{\jobname}

\def\cprime{$'$}
\begin{thebibliography}{10}

\bibitem{BlochOk}
S.~Bloch and A.~Okounkov.
\newblock The character of the infinite wedge representation.
\newblock {\em Adv. Math.}, 149(1):1--60, 2000.

\bibitem{Brown:thesis}
F.~Brown.
\newblock Multiple zeta values and periods of moduli spaces
  {$\overline{\mathscr M}_{0,n}$}.
\newblock {\em Ann. Sci. \'Ec. Norm. Sup\'er. (4)}, 42(3):371--489, 2009.

\bibitem{Brown:MMV}
F.~Brown.
\newblock Multiple modular values and the relative completion of the
  fundamental group of $\mathcal{M}_{1,1}$.
\newblock arXiv:1407.5167v3, 2016.

\bibitem{Brown:EquivIterIntEis}
F.~Brown.
\newblock A class of non-holomorphic modular forms {II}: Equivariant iterated
  eisenstein integrals.
\newblock arXiv:1708.03354, 2017.

\bibitem{DDMS}
M.~Deneufch{\^a}tel, G.~H.~E. Duchamp, V.~H.~N. Minh, and A.~I. Solomon.
\newblock Independence of hyperlogarithms over function fields via algebraic
  combinatorics.
\newblock In {\em Algebraic informatics}, volume 6742 of {\em Lecture Notes in
  Comput. Sci.}, pages 127--139. Springer, Heidelberg, 2011.

\bibitem{Eichler}
M.~Eichler.
\newblock Eine {V}erallgemeinerung der {A}belschen {I}ntegrale.
\newblock {\em Math. Z.}, 67:267--298, 1957.

\bibitem{Hai}
R.~Hain.
\newblock The geometry of the mixed {H}odge structure on the fundamental group.
\newblock In {\em Algebraic geometry, {B}owdoin, 1985 ({B}runswick, {M}aine,
  1985)}, volume~46 of {\em Proc. Sympos. Pure Math.}, pages 247--282. Amer.
  Math. Soc., Providence, RI, 1987.

\bibitem{Hain:Modulispaces}
R.~Hain.
\newblock Lectures on moduli spaces of elliptic curves.
\newblock In {\em Transformation groups and moduli spaces of curves}, volume~16
  of {\em Adv. Lect. Math. (ALM)}, pages 95--166. Int. Press, Somerville, MA,
  2011.

\bibitem{Hain:HodgeDeRham}
R.~Hain.
\newblock The {H}odge--de {R}ham theory of modular groups.
\newblock In {\em Recent advances in {H}odge theory}, volume 427 of {\em London
  Math. Soc. Lecture Note Ser.}, pages 422--514. Cambridge Univ. Press,
  Cambridge, 2016.

\bibitem{KanekoZagier}
M.~Kaneko and D.~Zagier.
\newblock A generalized {J}acobi theta function and quasimodular forms.
\newblock In {\em The moduli space of curves ({T}exel {I}sland, 1994)}, volume
  129 of {\em Progr. Math.}, pages 165--172. Birkh\"auser Boston, Boston, MA,
  1995.

\bibitem{Lang}
S.~Lang.
\newblock {\em Introduction to modular forms}.
\newblock Springer-Verlag, Berlin-New York, 1976.
\newblock Grundlehren der mathematischen Wissenschaften, No. 222.

\bibitem{LMS}
P.~Lochak, N.~Matthes, and L.~Schneps.
\newblock {Elliptic multiple zeta values and the elliptic double shuffle
  relations}.
\newblock arXiv:1703.09410, 2017.

\bibitem{Man}
Y.~I. Manin.
\newblock Iterated integrals of modular forms and noncommutative modular
  symbols.
\newblock In {\em Algebraic geometry and number theory}, volume 253 of {\em
  Progr. Math.}, pages 565--597. Birkh\"auser Boston, Boston, MA, 2006.

\bibitem{MartinRoyer}
F.~Martin and E.~Royer.
\newblock Formes modulaires et p\'eriodes.
\newblock In {\em Formes modulaires et transcendance}, volume~12 of {\em
  S\'emin. Congr.}, pages 1--117. Soc. Math. France, Paris, 2005.

\bibitem{MM}
J.~W. Milnor and J.~C. Moore.
\newblock On the structure of {H}opf algebras.
\newblock {\em Ann. of Math. (2)}, 81:211--264, 1965.

\bibitem{Radford}
D.~E. Radford.
\newblock A natural ring basis for the shuffle algebra and an application to
  group schemes.
\newblock {\em J. Algebra}, 58(2):432--454, 1979.

\bibitem{Reu}
C.~Reutenauer.
\newblock {\em Free {L}ie algebras}, volume~7 of {\em London Mathematical
  Society Monographs. New Series}.
\newblock The Clarendon Press, Oxford University Press, New York, 1993.
\newblock Oxford Science Publications.

\bibitem{Weibel:Hom}
C.~A. Weibel.
\newblock {\em An introduction to homological algebra}, volume~38 of {\em
  Cambridge Studies in Advanced Mathematics}.
\newblock Cambridge University Press, Cambridge, 1994.

\bibitem{123}
D.~Zagier.
\newblock Elliptic modular forms and their applications.
\newblock In {\em The 1-2-3 of modular forms}, Universitext, pages 1--103.
  Springer, Berlin, 2008.

\end{thebibliography}
\end{document}